\documentclass{article}
\usepackage{arxiv_style,times}
\usepackage{natbib}
\usepackage[utf8]{inputenc} 
\usepackage[T1]{fontenc}    
\usepackage{hyperref}
\hypersetup{colorlinks=true,linkcolor=blue,citecolor=green,urlcolor=blue}
\usepackage{url}            
\usepackage{booktabs}       
\usepackage{amsfonts}       
\usepackage{nicefrac}       
\usepackage{microtype}      
\usepackage{xcolor}         

\usepackage{url, subfiles, times}
\usepackage{dsfont}
\usepackage{xspace}
\usepackage{float,pgfplots,wrapfig,sidecap,lipsum}
\usepackage{booktabs,paralist}
\usepackage{amsfonts,  amsmath, amssymb, csquotes,xcolor,enumitem}
\usepackage{mathtools}
\usepackage[title]{appendix}
\usepackage{tablefootnote}
\usepackage{hyperref,psfrag}
\usepackage{tikz}
\usetikzlibrary{fit}
\usetikzlibrary{calc,shapes}
\usetikzlibrary{decorations.pathmorphing} 
\usetikzlibrary{fit} 
\usetikzlibrary{backgrounds}	
\usepackage{color,bm,mathrsfs,stmaryrd}
\usepackage[compact]{titlesec}
\usepackage{amsthm}
\usepackage{subcaption}
\usepackage{caption}

\newcommand{\F}{\mathscr{F}}

\newcommand{\B}{\mathbb{B}}
\newcommand\norm[1]{\left\lVert#1\right\rVert}
\newcommand\R{\mathbb{R}}
\newcommand\Z{\mathbb{Z}}

\newtheorem{theorem}{Theorem}
\newtheorem{lemma}[theorem]{Lemma}
\newtheorem{definition}{Definition}

\newtheorem{corollary}[theorem]{Corollary}

\newtheorem{assumption}{Assumption}

\usetikzlibrary{decorations.pathreplacing,calc}
\newcommand{\tikzmark}[1]{\tikz[overlay,remember picture] \node (#1) {};}
\newcommand*{\AddNote}[4]{%
    \begin{tikzpicture}[overlay, remember picture]
        \draw [decoration={brace,amplitude=0.5em},decorate,ultra thick,black]
            ($(#3)!(#1.north)!($(#3)-(0,1)$)$) --  
            ($(#3)!(#2.south)!($(#3)-(0,1)$)$)
                node [align=center, text width=1.85cm, pos=0.5, anchor=west] {#4};
    \end{tikzpicture}
}

\floatname{algorithm}{Algorithm}

\newcommand{\thm}[1]{\hyperref[thm:#1]{Theorem~\ref*{thm:#1}}}
\newcommand{\lem}[1]{\hyperref[lem:#1]{Lemma~\ref*{lem:#1}}}

\newtheorem{innercustomthm}{Theorem}
\newenvironment{customthm}[1]{\renewcommand\theinnercustomthm{#1}\innercustomthm}
{\endinnercustomthm}

\newtheorem{innercustomlem}{Lemma}
\newenvironment{customlem}[1]{\renewcommand\theinnercustomlem{#1}\innercustomlem}
{\endinnercustomlem}

\titlespacing{\section}{0pt}{1ex}{1ex}
    \titlespacing{\subsection}{0pt}{1ex}{0ex}
    \titlespacing{\subsubsection}{0pt}{0.5ex}{0ex}
    
\titlespacing{\paragraph}{%
  0pt}{
  0pt}{
  0.5em}
\pgfplotsset{compat=1.17}

\RequirePackage{algorithm}
\RequirePackage{algorithmic}

\def\mytitle{Escaping From Saddle Points Using \\ Asynchronous Coordinate Gradient Descent}
\def\mytitlerunning{Escaping From Saddle Points Using Asynchronous Coordinate Gradient Descent}

\title{\mytitle}

\author{Marco Bornstein, Jin-Peng Liu, Jingling Li, Furong Huang \\
Department of Computer Science, University of Maryland\\
\texttt{\{marcob, jliu1219, jingling, furongh\}@umd.edu} \\
}

\iclrfinalcopy 

\begin{document}

\maketitle

\begin{abstract}
   Large-scale non-convex optimization problems are expensive to solve due to computational and memory costs. 
 To reduce the costs, first-order (computationally efficient) and asynchronous-parallel (memory efficient) algorithms are necessary to minimize non-convex functions in machine learning. 
 However, asynchronous-first-order methods applied within non-convex settings run into two difficulties: (i) parallelization delays, which affect convergence by disrupting the monotonicity of first-order methods, and (ii) sub-optimal saddle points where the gradient is zero.
 To solve these two difficulties, we propose an asynchronous-coordinate-gradient-descent algorithm shown to converge to local minima with a bounded delay.
 Our algorithm overcomes parallelization-delay issues by using a carefully constructed Hamiltonian function.
 We prove that our designed kinetic-energy term, incorporated within the Hamiltonian, allows our algorithm to decrease monotonically per iteration.
 Next, our algorithm steers iterates clear of saddle points by utilizing a perturbation sub-routine.
 Similar to other state-of-the-art (SOTA) algorithms, we achieve a poly-logarithmic convergence rate with respect to dimension.
 Unlike other SOTA algorithms, which are synchronous, our work is the first to study how parallelization delays affect the convergence rate of asynchronous first-order algorithms. 
 We prove that our algorithm outperforms synchronous counterparts under large parallelization delays, with convergence depending sublinearly with respect to delays.
 To our knowledge, this is the first local optima convergence result of a first-order asynchronous algorithm for non-convex settings.
 
\end{abstract}

\setlength\abovedisplayskip{0pt}
\setlength\belowdisplayskip{0pt}

\section{Introduction}
\label{sec:intro}
First-order gradient based methods are widely used in large-scale optimization problems due to their computational efficiency compared with higher-order methods, such as Hessian-based methods. 
Recent advances have shown that these methods are also suitable for use in non-convex settings, with ~\citep{jin2017escape,jin2017accelerated, jin2019nonconvex} providing efficient local optima convergence guarantees. 
Even with these efficient first-order methods, optimization for ultra-large-scale non-convex problems, over-parameterized deep neural networks for instance ~\citep{Szegedy_2015_CVPR, zagoruyko2017wide}, is still computationally restrictive when solving on a single work node.
Thus, first-order methods commonly run in parallel across multiple workers to efficiently solve such large-scale problems.

To the best of our knowledge, only synchronous first-order methods have proven local optima convergence results within non-convex settings.
The requirement of synchronization, however, is impractical.
While the parallel computing setting allows for a speed-up in computational efficiency, it also comes with the potential dangers of parallelization delays and subsequent convergence slow-downs. 
The speed of parallel computing is limited by the weakest link in the computational chain: the slowest worker and its consequent longest communication delay. 
Worse, if a worker stops working or experiences a network connection failure, then the parallel computing process pauses.
In many applications, an asynchronous process can be implemented to train large-scale problems more efficiently than its synchronous counterpart \citep{chang2016asynchronous1,chang2016asynchronous2, assran2020advances}. 

Within non-convex settings, convergence to local optima is sufficient for many applications.
For non-convex problems, such as dictionary learning \citep{sun_dictionary}, matrix sensing and completion \citep{Bhoj2016, ge2016, park17a}, and tensor decomposition \citep{ge2015escaping}, all local optima are close in value to the global minimum. 
Therefore, convergence to local optima is a sufficient solution for these problems \citep{choromanska15}.
However, convergence to saddle points, or first-order stationary points in general, is insufficient.
Saddle points may produce sub-optimal solutions to a given optimization problem, such as training a deep neural network \citep{dauphin2014}. 

Asynchronous-first-order methods applied within non-convex settings run into two difficulties: (i) parallelization delays, which slow down convergence \citep{hannah2018unbounded, peng2019convergence}, and (ii) sub-optimal saddle points where the gradient is zero. 
In this paper, we present a first-order asynchronous coordinate gradient descent (ACGD) algorithm that overcomes these issues, and efficiently minimizes a high-dimensional non-convex objective function in the parallel computing setting.
Our algorithm overcomes delay issues by reformulating a non-convex function into a Hamiltonian, restoring the monotonicity of gradient methods lost in the presence of parallelization delays.
Furthermore, our algorithm avoids saddle points in a modified version of ~\citep{jin2017escape}: adding a carefully generated perturbation to the current iterate when the Hamiltonian does not sufficiently decrease to escape from potential saddle points. 
We prove later on that this generated perturbation, coupled with subsequent updates of the algorithm, dislodges iterates from a saddle point with high probability and sufficiently decreases the Hamiltonian.

\textbf{Summary of Contributions.}
To the best of our knowledge, this is the first local optima convergence result regarding a first-order asynchronous gradient based algorithm for non-convex settings with bounded parallelization delays.
Our main technical contributions are summarized below:\\
\textbf{(1)} Propose the first ACGD algorithm that converges $\epsilon$-close to second-order stationary points (properly defined in Section \ref{sec:preliminaries}) with high probability;\\
\textbf{(2)} Construct a novel Hamiltonian to restore monotonicity of ACGD under bounded delays;\\
\textbf{(3)} Determine how convergence is affected by bounded delays (in a sub-linear manner);\\
\textbf{(4)} Theoretically prove that the convergence of our algorithm is faster than synchronous counterparts under large parallelization delays.
\section{Related Work}
\label{sec:relatedwork}

\textbf{Convergence of First-order Algorithms.} Previously, it was believed that only second-order optimization methods would converge to second-order stationary points (such as cubic regularization or trust region algorithms). These algorithms are typically more computationally expensive \citep{gratton2016, nesterov2006} than recent, and more efficient, lower-order methods capable of converging to second-order stationary points \citep{ge2015escaping, Levy16a}. In these works, stochastic gradient descent is shown to converge to second-order stationary points, but with a high-degree polynomial dependence on the dimension. This result was improved in the work of Jin et al. \citep{jin2017escape,jin2017accelerated,jin2019nonconvex}. In \citep{jin2017escape,jin2019nonconvex}, perturbed stochastic and regular gradient descent converge to second-order stationary points poly-logarithmically with respect to dimension. Impressively, this mirrors the convergence rate to first-order stationary points (disregarding logarithmic factors). Our work matches this result, also relying poly-logarithmically on dimension.


Other simple first-order procedures that efficiently find an escaping direction from saddle points include NEON, within ~\cite{xu2017first} and ~\cite{allen2017neon2}. The work of ~\cite{xu2017first} is inspired by the perturbed gradient method proposed in~\cite{jin2017escape} and its connection with the power method for computing the largest eigenvector of a matrix starting from a random noise vector. In ~\cite{allen2017neon2}, negative-curvature-search subroutines are converted into first-order processes.
Our algorithm is complementary of this line of work, as one could combine these careful estimations of negative curvature into our algorithm.
The focus of our paper, rather than designing sophisticated procedures to improve the rate of escaping from saddle points, is to provide a simple algorithm that is resilient to delayed gradients.

\textbf{Convergence of Coordinate-Gradient-Descent Algorithms.} 
First mentioned by \citep{powell1973search}, coordinate gradient descent (CGD) has convergence difficulties in non-convex settings. Powell's finding that cyclic CGD fails to converge to a stationary point illustrates that a general convergence result for non-convex functions cannot be expected. However, Powell also mentions that the cyclic behavior of CGD is unstable with respect to small perturbations. To solve this issue, our paper follows \citep{sun2017asynchronous} in using inexact line searches to lead to convergence of CGD. Other methods and assumptions for attaining convergence results are discussed in \citep{wright2015}. These include assuming unique minimizers along any coordinate direction \citep{bertsekas1997nonlinear} and using functions with Kurdyka-Łojasiewicz (KL) properties \citep{attouch2010proximal,xu2017globally}.

\textbf{Convergence of Asynchronous Algorithms.} This paper follows the theory behind asynchronous coordinate gradient descent. A major portion of the theory, on asynchronous algorithms, has been built from \citep{arjevani20a,cole2018analysis, liu2015asynchronous,liu2014asynchronous}. The relationship between iterate delays and the convergence of asynchronous gradient methods is examined in \citep{arjevani20a}. It is shown in \citep{arjevani20a} that for convex quadratic functions, standard gradient descent methods have at most a linear dependence on delays. However, stochastic delayed gradient descent is found to have a smaller dependence on the delay. In \citep{liu2014asynchronous}, the convergence properties of an asynchronous-stochastic-coordinate-descent algorithm is nailed down for a convex function. Similar work is applied in \citep{liu2015asynchronous}, with convergence results for an asynchronous-stochastic-proximal-coordinate-descent algorithm. This algorithm has an improved complexity compared to similar asynchronous algorithms. Convergence analysis of asynchronous block coordinate descent in \citep{liu2015asynchronous,liu2014asynchronous} relied upon the independence assumption as well as bounded delays. Our work also incorporates the use of bounded delays, however without the need of the independence assumption. The work in \citep{cole2018analysis} uses a similar asynchronous, but accelerated, stochastic-coordinate-descent algorithm. While we do not incorporate acceleration, our work provides an even further speed up in convergence for convex functions than in \citep{liu2015asynchronous,liu2014asynchronous}.

The work of \citep{cole2018analysis, liu2015asynchronous,liu2014asynchronous} led to expanding research about asynchronous-coordinate-descent convergence within non-convex settings. In more closely related research to ours, the convergence properties of asynchronous coordinate descent in non-convex settings has been studied in \citep{sun2017asynchronous, cannelli2016asynchronous}. Within the important work of \citep{sun2017asynchronous}, asynchronous-coordinate-descent convergence results are determined for bounded, stochastic unbounded, and deterministic unbounded delays. These results are provided for both convex and non-convex functions. Furthermore, within \citep{sun2017asynchronous} it is shown that there is sufficient descent for bounded delays by using a Lyapunov function for a non-convex function. As stated above, \citep{sun2017asynchronous} proved that asynchronous coordinate gradient descent converges to stationary points for $L$-gradient Lipschitz non-convex functions (defined below). Our work utilizes the sufficient descent and Lyapunov function improvements found in \citep{sun2017asynchronous}. This paper dives deeper, and determines the relationship between bounded delays and convergence.
\section{Problem Formulation}
\label{sec:preliminaries}
\textbf{Notation.} We use $||\cdot||_2$ to denote the $\ell_2$-norm for vectors. 
We use $\lambda_{min}(\cdot)$ to denote the smallest eigenvalue of a matrix. 
For a function $f: \R^d \rightarrow \R$, the terms $\nabla f(\cdot)$ and $\nabla^2 f(\cdot)$ denote its gradient and Hessian respectively.
We denote the global minimum of $f$ as $f^*$, and $\Delta_f$ to be the difference between the function value at an initial solution $x^0$ and the global minimum, i.e., $\Delta_f := f(x^0) - f^*$. 
We use the notation $\mathcal{O}(\cdot)$ to hide absolute constants, and notation $\tilde{\mathcal{O}}(\cdot)$ to hide absolute constants and log factors.
Finally, we let $\B_x(r)$ denote the $d$-dimensional ball centered at $x$ with radius $r$.

\textbf{Problem Setting.} We solve an unconstrained, high dimension $d$, non-convex optimization
\begin{equation}
\setlength\abovedisplayskip{0pt}
\setlength\belowdisplayskip{0pt}
\label{eq:general_problem}
\min_{{x} \in \R^d} f({x}).
\end{equation}

\begin{assumption}
\label{assump:strict_saddle_lipschitz}
The objective function $f({x})$ is $L$-gradient Lipschitz, $\rho$-Hessian Lipschitz, $(\phi, \gamma, \zeta)$-strict-saddle, and has a bounded minimum value $f^*$. 
\end{assumption}
\begin{definition} A differentiable function $f(\cdot)$ is \textbf{$L$-smooth (or $L$-gradient Lipschitz)} if:
$$\forall x_1, x_2, \norm{\nabla f(x_1) - \nabla f(x_2)} \leq L \norm{x_1 - x_2}$$
\end{definition}

\begin{definition} A twice-differentiable function $f(\cdot)$ is \textbf{$\rho$-Hessian Lipschitz} if: 
$$\forall x_1, x_2, \norm{ \nabla^2 f(x_1) - \nabla^2 f(x_2) } \leq \rho \norm{x_1 - x_2} $$
\end{definition}

\begin{definition} \label{def:strict_saddle_function} A function $f(\cdot)$ is \textbf{$(\phi, \gamma, \epsilon)$-strict-saddle} if, for any $x$, at least one of the following holds:
$$||\nabla f(x)|| \geq \phi, \quad \lambda_{min}(\nabla^2 f(x)) \leq - \gamma, \quad x \text{ is } \epsilon\text{-close to } x^* - \text{ the set of local minima}$$
\end{definition}

\begin{definition}
\label{def:epsilon_second_order}
For a $\rho$-Hessian Lipschitz function $f(\cdot)$, we say that $x$ is a \textbf{second-order stationary point} if:
$$\norm{\nabla f(x)} = 0, \text{ and } \lambda_{min}(\nabla^2 f(x)) \geq 0$$ 
We also say x is an \textbf{$\epsilon$-second-order stationary point} if:
$$ \norm{\nabla f(x)} \leq \epsilon, \text{ and } \lambda_{min}(\nabla^2 f(x)) \geq -\sqrt{\rho \epsilon} $$
\end{definition}

\begin{definition} \label{def:strict_saddle} One calls $x$ a \textbf{saddle point} if $x$ is a first-order stationary point but not a local minimum. For a twice-differentiable function $f(\cdot)$, we say a saddle point $x$ is \textbf{strict (or non-degenerate)} if:
$$\lambda_{min}(\nabla^2 f(x)) < 0$$
\end{definition}
In many applications, such as dictionary learning, tensor decomposition, and matrix sensing, the objective function $f$ induces saddle points that are all strict \citep{sun_dictionary, Bhoj2016, ge2016, ge2015escaping}.
 
\textbf{ACGD Setting.} We seek to solve Equation \eqref{eq:general_problem} in an asynchronous coordinate gradient descent (ACGD) manner which efficiently solves high-dimensional problems by utilizing parallel computing without the costly synchronization bottlenecks.
ACGD differs from synchronous coordinate gradient descent (which is equivalent to parallel gradient descent), and regular gradient descent. 
While synchronous coordinate gradient descent can obtain a speed-up in convergence compared to the simpler regular gradient descent, its setting is more complex and comes with added difficulties: communication involved between the parameter server and workers within the parallel setting causes potential parallelization delays and subsequent convergence slow-downs. 
ACGD seeks to solve these issues that arise in synchronous coordinate gradient descent, and thus takes into account parallelization delays found within the parallel setting. 
Therefore, ACGD tackles the most difficult setting of the three: parallel gradient descent in the presence of parallelization delays.

Within the ACGD setting, there are $W$ workers, managed by a central parameter server, solving Equation~\eqref{eq:general_problem} in parallel without synchronization.
The parameter server assigns each worker $i$ a coordinate block $b_i$ to update, where $\bigcup\limits_{i=1} b _i = \{x_1,x_2,\ldots,x_d\} = x \in \R^d$. The ACGD process proceeds as follows:\\
\textbf{(1)} Worker $i$ fetches the current up-to-date global iterate $x^j$ from the parameter server, where $j$ indicates the current global iteration number that the parameter server monitors. \\
\textbf{(2)} Worker $i$ then begins computing the gradient $\nabla_{b_i}f(x^{j})$  corresponding to its coordinate block $b_i$.\\
\textbf{(3)} Once worker $i$ finishes, it sends its block of gradient values $\nabla_{b_i}f(x^{j})$ back to the server. \\
\textbf{(4)} The server receives the block gradient from worker $i$, yet the global iteration $j$ has been \emph{updated} since worker $i$ previously fetched the most up-to-date iterate in \textbf{(1)}. As other workers have been updating the global iterate simultaneously, the once up-to-date iterate that worker $i$ received in \textbf{(1)} has now decayed to become out-of-date after completion of \textbf{(3)}. To reflect this, the block gradient the server receives from worker $i$ is denoted as $\nabla_{b_i}f(\hat{x}^{j}_i)$, since $j$ has been updated to a different value (where $\hat{x}^{j}_i$ indicates the out-of-date iterate for worker $i$ as will be defined in Equation~\eqref{eq:delay}). The parameter server uses $\nabla_{b_i}f(\hat{x}^{j}_i)$ to update the coordinates of the global iterate corresponding to $b_i$ with a step size $\eta$:
\begin{equation}\label{eq:ACGD}
    x^j_{b_i} \gets x^j_{b_i} - \eta \nabla_{b_i}f(\hat{x}^{j}_i)
\end{equation} and the global iteration number increments $j \gets j +1$. \\
\textbf{(5)} The entire updated global iterate $x^j$ is sent back just to worker $i$ so it can repeat steps \textbf{(2)-(4)}.

\begin{figure}[!htbp]
\centering
\begin{algorithm}[H]
\caption{SW-ACGD($x,f,\eta,b$): Single-Worker ACGD}
\label{alg:Algo1}
\begin{algorithmic}[1]
    \STATE {\bfseries Input:} current iterate $x$,  function $f$, step size $\eta$, assigned coordinate block $b$
    \STATE $u \gets \boldsymbol{0} = {0}^{d}$
    \STATE $u_{b} \gets -\eta \nabla_{b} f(x)$
    \STATE {\bfseries Output:} $u$
\end{algorithmic}
\end{algorithm}
\end{figure}


\textbf{Asynchronous Updates and Bounded Delays.} As mentioned in the ACGD process above, the global iteration number $j$ is updated when one of the workers $i \in \{1, \ldots, W \}$ sends back its block of gradients to the parameter server. 
In step \textbf{(4)} above, the global iteration number $j$ may have been updated, i.e., the values of $j$ in step \textbf{(1)} and step \textbf{(4)} may be different, since the parameter server may have received updates from other workers between the time worker $i$ fetches the current global iterate and the time worker $i$ sends back updates. 
Equation \eqref{eq:ACGD} uses $\hat{x}^{j}_i$ rather than $x^{j}$ when computing the gradient. 
This is where ACGD differs from gradient descent and synchronized block coordinate gradient descent. 
In the asynchronous parallel setting, workers are able to compute their block of gradients, send them to the parameter server, and receive an updated global iterate without waiting for other slower workers.
Therefore, slower workers will be left computing the gradient of iterates that are out of date compared to the global iterate. 
Delays arise from these slow workers or workers with expensive blocks to update (blocks with many non-zero entries compared to others). 
The term $\hat{x}^{j}_i$ denotes this out-of-date iterate for worker $i$, and is defined as:
\begin{equation}
\label{eq:delay}
    \hat{x}^{j}_i := \bigg(x_1^{j - D_i(1)}, \; \ldots \; x_i^{j} \; \ldots, \; x_W^{j - D_i(W)}\bigg).
\end{equation}

\begin{figure}
   \centering
    \includegraphics[width = 0.9\textwidth]{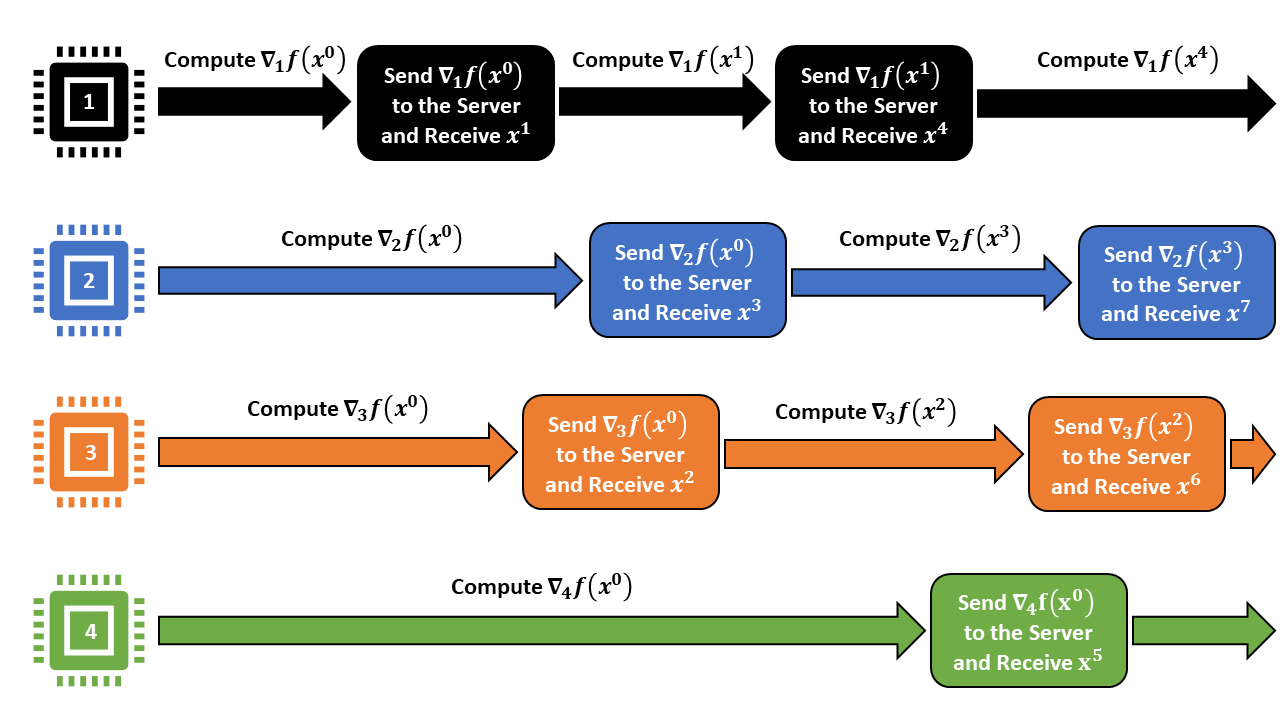}
    \vspace{-1em}
    \caption{ACGD Schematic with Four Workers.}
    \label{fig:ACGD-Process}
\end{figure}
As one can see above, $\hat{x}^{j}_i$ is split up into $W$ blocks, a block corresponding to each worker.
The term $D_i(p) \in \Z^{+}$ is defined as the number of iterations elapsed since worker $p$ communicated its \textit{first} coordinate update to the parameter server throughout the duration worker $i$ is performing its current update.
Thus, $D_i(p)$ denotes the delay at the $p^{\text{th}}$ coordinate block for worker $i$. 
The maximum bounded delay $\tau$, a user input, is an upper bound on all worker delays: $\max_{1\leq p \leq W} \{ D_i(p) \} \leq \tau, \; \forall \; i$. 
In Equation \eqref{eq:delay}, $\hat{x}_i^j$ has no delay term on the $i$-th block as only worker $i$ updates the $i$-th block.




It is important to note that ACGD experiences zero delay if $\tau = W-1$. 
First, it is impossible for $\tau$ to be smaller than $W-1$. 
With $W$ workers, updating each coordinate block of the solution would take a minimum of $W$ global iterations. 
In this perfect scenario, which only occurs if the blocks are updated in a cyclic manner with no randomness, the final coordinate block waiting to be updated would have waited $W-1$ iterations. 
Thus, the zero-delay scenario occurs when $\tau = W-1$.

The ACGD algorithm we propose applies to both zero-delayed ($\tau = W-1$) and delayed settings ($\tau > W-1$). 
Thus, we allow $\tau \geq W-1$.
In practice, and in our analysis, we set $\tau \geq 1$ as we assume multiple parallel workers exist. 
If $W=1$, the asynchronous parallel setting trivially reduces to the sequential gradient descent setting. 
Furthermore, if $W = 1$, our ACGD algorithm reduces to the perturbed gradient descent algorithm from \citep{jin2017escape}. 
\section{Hamiltonian: A Modified Objective Function for Monotonic Descent}
\label{sec:Hamiltonian}

Bounded delays arise from the setting of ACGD, as detailed in Section \ref{sec:preliminaries}, and adversely affect convergence analysis of ACGD algorithms. 
Namely, under parallelization delays, the objective function ceases to decrease monotonically. 
Even with a small enough step-size, ACGD is not guaranteed to be monotonically decreasing.
This leads to difficulty in convergence analysis.
One key technical contribution from our paper is proposing a Hamiltonian (sum of potential energy and kinetic energy) \citep{jin2017accelerated,sun2017asynchronous}, which has a monotonic property, to guarantee convergence of ACGD. 
 

\textbf{Construction of the Discrete Hamiltonian.} We show below that a discrete Hamiltonian can be created to ensure monotonicity of Equation \eqref{eq:ACGD}, the discretized version of ACGD. 
Similar to the work of \cite{sun2017asynchronous}, we utilize a Lyapunov function to model the discrete version of ACGD.
Let the bounded delays be $\tau \geq 1$, $\tau \in \Z^+$. We now define the discrete Hamiltonian as the following:
\begin{definition}[\textbf{Discrete Hamiltonian}]
\label{def:hamiltonian_discrete}
The Hamiltonian, with a bounded delay $\tau \geq 1$, is defined in the discrete case as 
\begin{equation}
    E(x^j) = E_j := f(x^j) + \frac{L}{2\sqrt{\tau}}\sum_{i=j - \tau}^{j-1}(i - (j - \tau) + 1)||x^{i+1} - x^i||_2^2.
\end{equation}
\end{definition}
As one can see above in Definition \ref{def:hamiltonian_discrete}, the Hamiltonian's potential energy is denoted as the true objective function value $f(x^j)$. 
Thus, when minimizing the Hamiltonian, we are still minimizing $f$.

\textbf{Guaranteeing Discrete Hamiltonian Monotonicity.} 
To guarantee monotonicity of the discrete Hamiltonian, it must be shown that each subsequent iterate of the Hamiltonian decreases in value. 
The first step in guaranteeing monotonicity is bounding the difference between consecutive iterates of the Hamiltonian below, and proving that this lower bound is greater than zero for all iterates. 
Lemma \ref{lem:bounded_delay_decrease} provides a result about this lower bound, and is crucial to our Main Theorem.
\begin{lemma}[Guaranteed Discrete Hamiltonian Decrease]
\label{lem:bounded_delay_decrease}
Under Assumption \ref{assump:strict_saddle_lipschitz}, let the sequence
$\{x^j\}_{j \geq 0}$ be generated by ACGD with delays bounded by $\tau \geq 1$. There exists $\beta > 0$ such that choosing the step size
    $\eta \leq \frac{1}{2L\tau^{1/2 - \beta} \iota \chi }$ 
results in
\begin{equation}
    E_j - E_{j+1} \geq L(\frac{1}{\eta L} - \sqrt{\tau} - \frac{1}{2})||x^{j+1} - x^{j}||_2^2 > 0.
\end{equation}
\end{lemma}
\textbf{Remark.} (1) The discrete Hamiltonian is guaranteed to monotonically decrease with a small step-size. (2) Similarly, the step-size is inversely proportional to the bounded delay, meaning descent slows with a large delay. (3) The discrete Hamiltonian decreases proportionally to the difference between consecutive iterates, and therefore converges to an optimal solution along with our iterates. (4) Step-size for sufficient descent mirrors the optimal step-size in analysis of gradient descent convergence to first-order stationary points $O(1/L)$ \citep{nesterov}.

The parameter $\beta$ determines the degree of how sub-linear the convergence of our proposed algorithm depends upon the maximum delay bound $\tau$. 
The hyperparameters $\iota \geq 3$ and $\chi \geq 1$ are defined below in Appendix \ref{appendix:parameters}. 
Calculating $\beta$ is shown in Equation \eqref{eq:beta} within Appendix \ref{appendix:parameters} (a simple optimization problem with a closed-form solution). For small values of $\tau$, $\beta$ assumes a value close to 1/2. 
Thus, for small bounded delay values, the step-size $\eta$ exactly matches up with \citep{nesterov} $O(1/L)$. 
As bounded delays grow, the value of $\beta$ decreases at a slow rate. Even for large value bounded delays ($\tau > 100$), $\beta$ remains large ($\beta >$ 1/9). 
Thus, our algorithm achieves a significant sub-linear convergence result regardless of the size of the bounded delays.

\begin{corollary}
\label{cor:suff_decrease}
With the step-size specified in Lemma \ref{lem:bounded_delay_decrease}, one finds: $E_j - E_{j+1} \geq \frac{3}{8}L||x^{j+1} - x^{j}||_2^2$.
\end{corollary}
\textbf{Remark.} For a small step size, the Hamiltonian decreases in proportion to only the difference in iterates and the gradient-Lipschitz, with no dependence on the bounded delays.
\section{Saddle-Escaping ACGD: Overview and Convergence Results}
\label{sec:Main_Result}

In this section, we present an overview of our algorithm, entitled Saddle-Escaping ACGD (SE-ACGD), and its efficient convergence to $\epsilon$-second-order stationary points (defined in Section \ref{sec:preliminaries}). 

\textbf{Saddle-Escaping ACGD Process.} The SE-ACGD algorithm works via the following steps:

\textbf{(1)} Begin with an initial iterate $x^0$ and Hamiltonian $E_0$.\\
\textbf{(2)} Perform the ACGD process for at most $\tau+1$ iterations, outputting $E_{lg}, x^{lg}$ (Algorithm \ref{alg:Algo2}).\\
\textbf{(3)} Repeat step (2) until the difference between the starting and ending Hamiltonians is small.\\
\textbf{(4)} Store the outputted Hamiltonian $E_s$ and iterate $x^s$ (an $\epsilon$-first-order stationary point) from step \textbf{(3)}.\\
\textbf{(5)} Perturb $x^s$ and run $T$ iterations of the ACGD process, outputting $E_{p}, x^{p}$ (Algorithm \ref{alg:Algo3}).\\
\textbf{(6)} Check if the Hamiltonian difference $E_s - E_{p}$ is small.\\
\textbf{(7)} If small, $x^{s}$ is returned (it is an $\epsilon$-second-order stationary point) and if not, steps \textbf{(2)-(6)} are repeated starting from $x^{p}$ (a saddle point is escaped).

We define the difference in Hamiltonians to be small if it is smaller than a specified function threshold $\F$. 
This threshold is a hyperparameter, defined in Appendix \ref{appendix:parameters}, which is our barometer to ensure that the Hamiltonian is decreasing towards the optimal value.
Below, we present the algorithmic layout of SE-ACGD.
Within SE-ACGD, two sub-routines are utilized.
The first, LG-ACGD, descends along large gradient regions (step \textbf{(2)}) and is introduced later in Algorithm \ref{alg:Algo2}. 
The second, P-ACGD, descends along saddle point regions (step \textbf{(5)}) and is introduced later in Algorithm \ref{alg:Algo3}.

\begin{figure}[!htbp]
\centering
\begin{algorithm}[H]
\caption{SE-ACGD($x^0, f, \tau, L, \eta, r, T, \F$):
Saddle-Escaping ACGD}
\label{alg:MainAlgo}
\begin{algorithmic}[1]
    \STATE {\bfseries Input:} initial iterate $x^0$, function $f$, bounded delay $\tau$, gradient Lipschitz $L$, $^*$step size $\eta$, $^*$perturbation radius $r$, $^*$perturbation steps $T$, $^*$Hamiltonian threshold $\F$
    \STATE $E_0 \leftarrow f(x^0), \; j \leftarrow 0$
    \FOR{$s = 1, 2, 3, \ldots $}
        \STATE $x^{j+\tau+1}, \; E_{j+\tau+1} \leftarrow$ LG-ACGD($x^j$, $j$, $f$, $\eta$, $\tau$, $L$)
        \STATE $j \leftarrow j + \tau + 1$
        \IF{$(E_{j-1} - E_j) < \F$}
            \STATE $x^{j+T}, E_{j+T} \leftarrow$ P-ACGD($x^j$, $f$, $\eta$, $\tau$, $r$, $T$, $L$)
            \IF{$(E_{j} - E_{j+T}) < \F$}
                \STATE STOP: $x^j$ is an $\epsilon$-Second-Order Stationary Point
            \ENDIF
            \STATE $j \leftarrow j + T$
        \ENDIF
    \ENDFOR
    \STATE {\bfseries Output:} $x^{j}$
    \STATE {\small $^*$ This hyperparameter is computed as in Appendix \ref{appendix:parameters}}
\end{algorithmic}
\end{algorithm}
\end{figure}

\textbf{Convergence Result and Comparison.} 
Below, we detail the convergence rate of SE-ACGD.
First, we mention how iteration complexities can vary in parallel settings. 
In parallel settings, updates from \textit{each} worker must be accounted for in the iteration complexity.
Due to bounded delays, convergence results of synchronous algorithms include a linear term $\tau$ since all workers wait for the slowest worker at each round.
Thus, the convergence result of parallel synchronous GD is $\tilde{\mathcal{O}}(\tau L\Delta_f/\epsilon^2)$.


\begin{theorem}[\textbf{Main Convergence Theorem}]\label{thm:main_theorem}
Let the function $f$ satisfy Assumption \ref{assump:strict_saddle_lipschitz}. For any $\delta>0$ and $\epsilon\le\frac{L^2}{\rho}$, there exists a step size $\eta$ satisfying Lemma \ref{lem:bounded_delay_decrease} such that Algorithm \ref{alg:MainAlgo} will output an $\epsilon$-second-order stationary point, with probability $1-\delta$, in the following number of iterations:
\begin{equation}
\label{eq:complexity}
    \mathcal{O}\bigg(\frac{L^2 \Delta_f \tau^{1-2\beta}}{\epsilon^{2.5}\rho^{0.5}} \log_2^4(\frac{d \tau L \Delta_f}{\delta \epsilon})\bigg).
\end{equation}
\end{theorem}

\textbf{Remark.} \textbf{(1)} 
Our ACGD algorithm obtains a faster convergence rate to $\epsilon$-second-order stationary points than existing parallel algorithms when $L/\tau^{2\beta}\sqrt{\rho \epsilon} \leq 1$. 
This occurs when $\tau$ is large, which happens when (i) there is a very slow worker and/or (ii) there are many workers (since $\tau \geq W-1$).
Our theoretical result makes sense intuitively: asynchronous algorithms are used in practice when (i) there are slow or unresponsive workers and/or (ii) numerous workers are needed to solve a large-scale problem.
\textbf{(2)} The convergence to $\epsilon$-second-order stationary points with an asynchronous process depends sub-linearly on the bounded delay $\tau$. 
This is an important first step in analyzing how bounded delays affect the iteration complexity for asynchronous algorithms in the non-convex setting. 

Our work is the first to prove the following intuition to be correct: asynchronous algorithms perform better than synchronous algorithms, under large delays, when converging to local minima in non-convex settings.
The sub-linear convergence dependence with respect to the maximum bounded delay is an improvement over the synchronous baseline. 
When delays are small (small $\tau$), they do not affect convergence.
This occurs, since small values of $\tau$ push $\beta$ to equal 1/2.
Finally, our algorithm does not match the rate on $\epsilon$ compared with~\cite{jin2017accelerated} since no acceleration is used.
This is left for future work.

\textbf{Proof Sketch.} To prove convergence to $\epsilon$-second-order stationary points, we prove that our SE-ACGD algorithm descends along regions with a large gradient and regions that are close to, or located at, a saddle point.
Once both regions have been descended, the set of local optima has been reached.
We prove that in large-gradient regions, where $||\nabla f(x^j)||_2 > \epsilon$, our algorithm sufficiently decreases the Hamiltonian, our modified objective function, by a sufficient amount $\F$ (a hyperparameter).
Likewise, we prove that in saddle-point regions, $||\nabla f(x^j)||_2 \leq \epsilon$ and $\lambda_{min}(\nabla^2 f(x^j)) < -\sqrt{\rho \epsilon}$, our algorithm decreases the Hamiltonian by $\F$ with a high probability.
Proving sufficient descent along saddle-point regions is complicated, necessitating a bound on the volume of the stuck region around saddle points (the space where points perturbed from a saddle point remain stuck at it after a certain number of ACGD iterations).
This is discussed in Section \ref{sec:Convergence_Analysis} and in full detail within Appendix \ref{appendix:Main-Theorem-Proof}.
\section{Saddle-Escaping ACGD: Algorithm and Attaining Convergence}
\label{sec:Convergence_Analysis}

In this section, we demonstrate how our SE-ACGD algorithm is constructed to ensure convergence to local optima, $\epsilon$-second-order stationary points, in non-convex settings. 

To begin, we remind the reader that the given objective function $f$ is assumed to be strict-saddle. If a function is strict-saddle, the landscape of the function can be split up into three regions. The first is a large-gradient region, the second is a region containing a notable negative eigenvalue of the Hessian (near a saddle point), and the third is a region close to a local minimum, or second-order stationary point. Our algorithm converges to a region $\epsilon$-close to a second-order stationary point by sufficiently descending along both large-gradient and saddle point regions. 

\textbf{Descending Along Large-Gradient Regions.} Descent along large-gradient regions, occurring in step \textbf{(2)} of the SE-ACGD process detailed in Section \ref{sec:Main_Result}, is accomplished by the sub-algorithm titled Large-Gradient ACGD (LG-ACGD), displayed in Algorithm \ref{alg:Algo2}. 

\begin{figure}[!htbp]
\centering
\begin{algorithm}[H]
\caption{LG-ACGD($x^j, j, f,\eta,\tau,L$):
Large-Gradient ACGD}
\label{alg:Algo2}
\begin{algorithmic}[1]
    \STATE {\bfseries Input:} current iterate $x^j$, global iteration $j$, function $f$, step size $\eta$, bounded delay $\tau$, gradient Lipschitz $L$
    \FOR{$s = 1, \ldots, \tau+1$\tikzmark{top3}}
        \STATE $x^{j+1} \gets x^j +$ SW-ACGD($x^j$, $f$, $\eta$, $b_i$)
        \STATE $j \leftarrow j + 1$\tikzmark{bottom3}
    \ENDFOR
    \STATE $E_{j} = f(x^j) + \frac{L}{2\sqrt{\tau}}\sum\limits_{k = j - \tau}^{j-1} (k - j + \tau +1)|\tikzmark{right3}|x^{k+1} - x^{k}||_2^2$
    \STATE {\bfseries Output:} $x^{j}$, $E_{j}$ 
    \AddNote{top3}{bottom3}{right3}{In Parallel\\For Any\\Block $b_i$} 
\end{algorithmic}
\end{algorithm}
\end{figure}

Within LG-ACGD, the ACGD process is run in parallel by calling SW-ACGD for each worker $i$.
This parallel process is cumulatively performed for $\tau + 1$ iterations.
The reason for $\tau + 1$ iterations, is that
all coordinate blocks will be updated within the window of $\tau + 1$ iterations due to the bounded delays assumption.
Once guaranteed that each coordinate block $b_i$ will be updated within $\tau + 1$ iterations, we can prove sufficient descent along large-gradient regions.

\begin{theorem}[\textbf{Large Gradient Scenario}]\label{thm:large_gradient}
Let function $f$ satisfy Assumption \ref{assump:strict_saddle_lipschitz}, $\eta$ satisfy Lemma \ref{lem:bounded_delay_decrease}, and let $\tau \geq 1$. If $||\nabla f(x^j)||_2 > \epsilon$, then by running Algorithm \ref{alg:Algo2} we have:
    $E_{j-\tau} - E_{j+1} > \F$.
\end{theorem}

The result of Theorem \ref{thm:large_gradient} solidifies the relationship between the sufficient decrease of the Hamiltonian and the bounded delays from the asynchronous process.
When the gradient is large, only after $\tau + 1$ updates from asynchronous workers will the Hamiltonian decrease by a specified threshold $\F$. 
Thus, Theorem \ref{thm:large_gradient} reiterates that for the Hamiltonian to sufficiently decrease, all asynchronous workers must finish computing their gradient blocks. The proofs of Theorem \ref{thm:large_gradient} is detailed in Appendix \ref{appendix:large_grad}.

\textbf{Descending Along Saddle-Point Regions.} Descent along saddle point regions, step \textbf{(5)} of the SE-ACGD process detailed in Section \ref{sec:Main_Result}, is performed by the sub-algorithm titled Perturbed ACGD (P-ACGD), shown in Algorithm \ref{alg:Algo3}.

\begin{figure}[!htbp]
\centering
\begin{algorithm}[H]
\caption{P-ACGD($x^j, f,\eta,\tau,r,T,L$):
Perturbed ACGD}
\label{alg:Algo3}
\begin{algorithmic}[1]
    \STATE {\bfseries Input:} current iterate $x^j$, function $f$, step size $\eta$, bounded delay $\tau$, perturbation radius $r$, perturbation steps $T$, gradient Lipschitz $L$
    \STATE $\xi \leftarrow \text{ uniformly} \sim \B_{x^j}(r)$
    \STATE $y^0 \leftarrow x^j + \xi$, $\; t \gets 0$ 
    \WHILE{$t < T$\tikzmark{top2}}
        \STATE $y^{t+1} \gets y^t +$ SW-ACGD($y^t$, $f$, $\eta$, $b_i$)
        \STATE $t \leftarrow t + 1$\tikzmark{bottom2}
    \ENDWHILE
    \STATE $E_{T} = f(y^T) + \frac{L}{2\sqrt{\tau}}\sum\limits_{k = T - \tau}^{T-1} (k - T + \tau + 1)\tikzmark{right2}||y^{k+1} - y^k||_2^2$
    \STATE $x^{j+1} \leftarrow y^T$, $E_{j+1} \leftarrow E_T$
    \STATE {\bfseries Output:} $x^{j+1}$, $E_{j+1}$
    \AddNote{top2}{bottom2}{right2}{In Parallel For Any Block $b_i$}
\end{algorithmic}
\end{algorithm}
\end{figure}

The P-ACGD algorithm is designed to escape saddle points. 
Called after multiple rounds of LG-ACGD, P-ACGD is utilized once an iterate $x^j$ is proven to be outside of the large-gradient region.
The algorithm begins by generating a perturbed iterate $y^0$, which is uniformly selected from a $d$-dimensional ball with radius $\eta r$ centered at $x^j$.
Finally, the ACGD process is run in parallel by calling SW-ACGD for each worker $i$.
This parallel process is cumulatively performed for $T$ (a hyperparameter) iterations.

It is proven below that after $T$ iterations, the carefully generated perturbed iterate will have a high probability of escaping a saddle point.
Therefore, the output of P-ACGD will be an iterate that escapes a saddle point, and decreases the Hamiltonian by $\F$, or an iterate that remains at a stationary point --- an $\epsilon$-second-order stationary point.

\begin{theorem}[\textbf{Saddle Point Scenario}]\label{thm:saddle_point}
Let function $f$ satisfy Assumption \ref{assump:strict_saddle_lipschitz} and let $\eta$ satisfy Lemma \ref{lem:bounded_delay_decrease}. If $\|\nabla f(x^j)\|_2 \leq \epsilon$ and $\lambda_{\min}(\nabla^2 f(x^j)) < -\sqrt{\rho\epsilon}$, then $x^j$ is located at a saddle point. Let $y^0=x^j+\xi$, where $\xi$ comes from the uniform distribution over ball with radius $\eta r$, and $\{y^t\}_{t=0}^T$ are the iterates of ACGD starting from $y^0$ with $T \geq \frac{\text{log}_2(\sigma \iota^2 \chi^2)}{\eta \sqrt{\rho \epsilon}}$. Let $x^{j+1} = y^T$. Then, with at least probability $1-\gamma$, running Algorithm \ref{alg:Algo3} once results in:
    $E_{j} - E_{j+1} > \F$.
\end{theorem}
\textbf{Proof Sketch.}  This proof is inspired by the geometry of saddle points \citep{jin2017escape,jin2017accelerated}. Let the point $x_s$ satisfy the conditions of Theorem \ref{thm:saddle_point} (a saddle point). To escape the saddle, one perturbs $x_s$. Adding a perturbation ($x_p = x_s + \xi$) results in $x_p$ being generated from a uniform distribution over a perturbation ball. This perturbation ball is centered at $x_s$ with a radius of $\eta r$. 

There are two regions within this perturbation ball: a stuck region and an escaping region. 
The stuck region refers to the region within the perturbation ball where running $T$ iterations of ACGD on $x_p$ will lead back to, and remain stuck at, the saddle point. Conversely, the escaping region refers to the region where running $T$ iterations of ACGD on $x_p$ will escape the saddle point. 

The stuck region is quite small compared to the rest of the perturbation ball \citep{jin2017escape,jin2017accelerated}. In fact, the width of the stuck region in the direction of the minimum eigenvalue is proven to be bounded by a small constant multiplied by the radius of the perturbation ball. This width is used to bound the volume of the stuck region, and we prove with high probability that $x_p$ is unlikely to fall within the stuck region. Our work builds from this geometric technique, and the full proof of Theorem \ref{thm:saddle_point} is in the appendix.

\textbf{Outputting An \texorpdfstring{$\epsilon$}{TEXT}-Second-Order Stationary Point.}
To summarize, our SE-ACGD algorithm works as follows. 
First, LG-ACGD is run on the initial iterate $x^0$ until the outputted iterate $x^s$ is no longer within the large-gradient region (the Hamiltonian does not sufficiently decrease).
Once this occurs, P-ACGD is run on $x^s$, where $x^s$ is perturbed and undergoes $T$ iterations of ACGD, producing an output $x^{p}$.
If, after $T$ iterations, the subsequent Hamiltonian sufficiently decreases, then $x^s$ was previously contained in the saddle-point region.
By Theorem \ref{thm:saddle_point}, the output $x^{p}$ of P-ACGD is not contained in the saddle-point region (and thus has escaped the saddle point located at $x^s$).
However, if, after $T$ iterations, the subsequent Hamiltonian does not sufficiently decrease, then $x^s$ is not contained in the saddle-point region.
Thus, $x^s$ is not located in either the large-gradient or the saddle-point regions. 
Since the objective function is strict-saddle, it follows that $x^s$ is an $\epsilon$-second-order stationary point.
Thus, SE-ACGD finds a local optimum.

\begin{figure}[!htbp]
    \centering
    \includegraphics[width=0.7\textwidth]{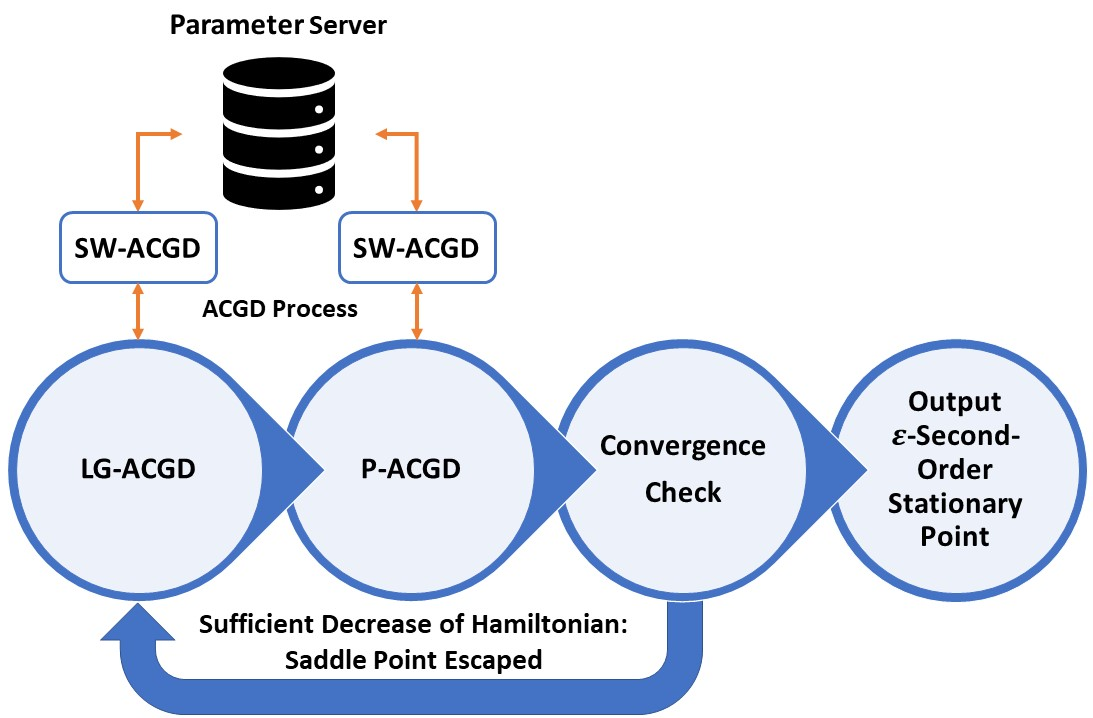}
    \caption{\textit{Saddle-Escaping ACGD Schematic.} Starting on the left, SE-ACGD calls LG-ACGD to find a stationary point $x^s$. Once found, P-ACGD is called to perturb and dislodge $x^s$ in case it is a saddle point. These sub-routines are powered by parallel workers performing ACGD and communicating with the server (SW-ACGD). If $x^s$ is a saddle point, the output $x^p$ of P-ACGD sufficiently reduces the Hamiltonian and the process restarts. If not, $x^s$ is an $\epsilon$-second-order stationary point.}
    \label{fig:seacgd}
\end{figure}

\section{Experiments}
\label{sec:experiments}
Below, we test SE-ACGD against other baseline first-order algorithms (serial gradient descent and parallel perturbed gradient descent) to analyze its performance in the following criteria:
\begin{enumerate}
    \item Effectiveness at escaping saddle points.
    \item Performance under parallelization delays.
    \item Parallel speed-up and scalability.
\end{enumerate}
The problem we seek to minimize is the following function:
\begin{align}
    f(\bm{x}) & := \tilde{f}(\frac{2}{d}\sum_{i=1}^{d/2} x_i, \frac{2}{d}\sum_{j=d/2+1}^d x_j), \quad \tilde{f}(r,s) := d\big[(r-1)^4 - (r-1)^2 + (s+1)^2\big].
\end{align}

The equation above has a saddle point located at $(r,s) = (1, -1)$ with a function value of 0.
Furthermore, there exist two equal local optima located at $(r,s) = (1 + \frac{1}{\sqrt{2}}, -1) = (1 - \frac{1}{\sqrt{2}}, -1)$.
These local optima have a minimum function value of $-d/4$. 
Below, we showcase how SE-ACGD reaches the local optima quick and efficiently.

\subsection{Saddle Point Evasion}
\label{app:saddle-evasion-exp}

Our first additional experiment analyzes how well SE-ACGD escapes saddle points compared to its peer algorithms.
Furthermore, we investigate whether SE-ACGD is effective at escaping saddle points in high dimensions.

\begin{figure}[!htbp]
     \centering
     \begin{subfigure}[b]{0.475\textwidth}
         \centering
         \includegraphics[width=0.75\textwidth]{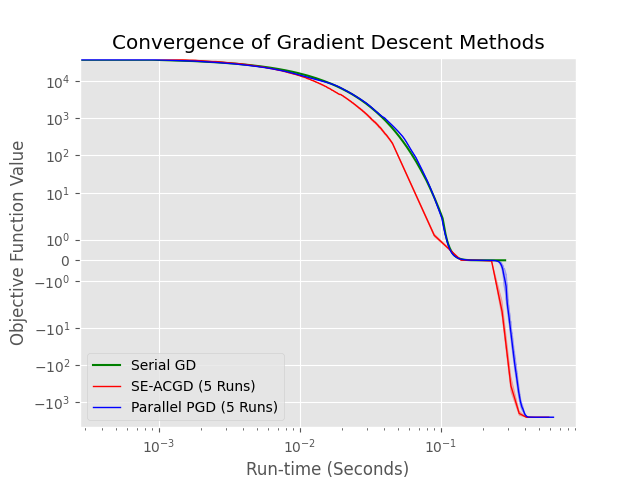}
         \caption{1e4 Dimensions, 8 Workers, No Delay.}
         \label{fig:8w-10k}
     \end{subfigure}
     \hfill
     \begin{subfigure}[b]{0.475\textwidth}
         \centering
         \includegraphics[width=0.75\textwidth]{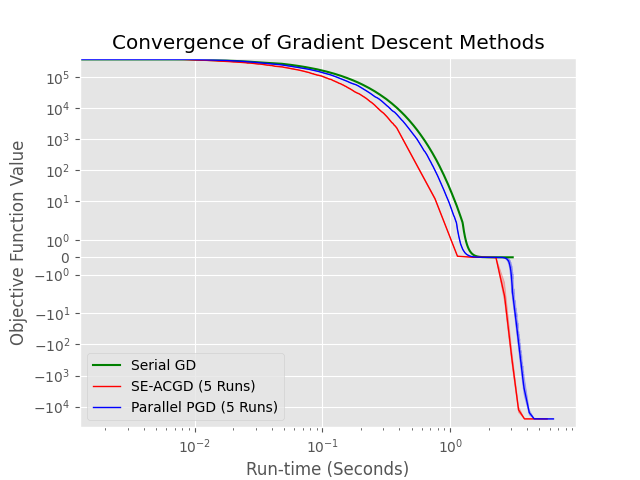}
         \caption{1e5 Dimensions, 8 Workers, No Delay.}
         \label{fig:fig:8w-100k}
     \end{subfigure}
     \hfill
     \begin{subfigure}[b]{0.475\textwidth}
         \centering
         \includegraphics[width=0.75\textwidth]{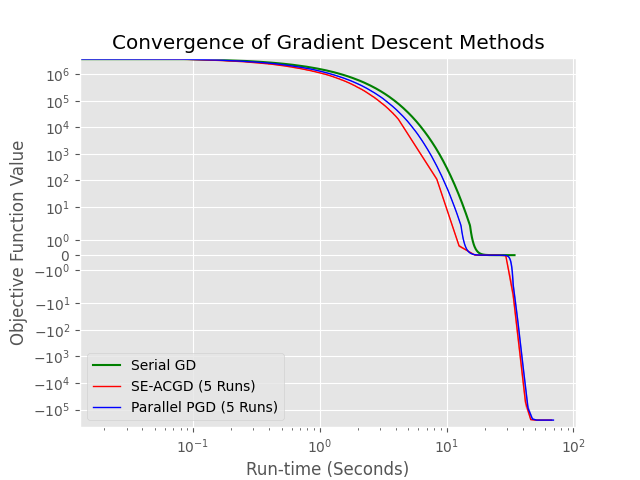}
         \caption{1e6 Dimensions, 8 Workers, No Delay.}
         \label{fig:8w-1M}
     \end{subfigure}
     \hfill
     \begin{subfigure}[b]{0.475\textwidth}
         \centering
         \includegraphics[width=0.75\textwidth]{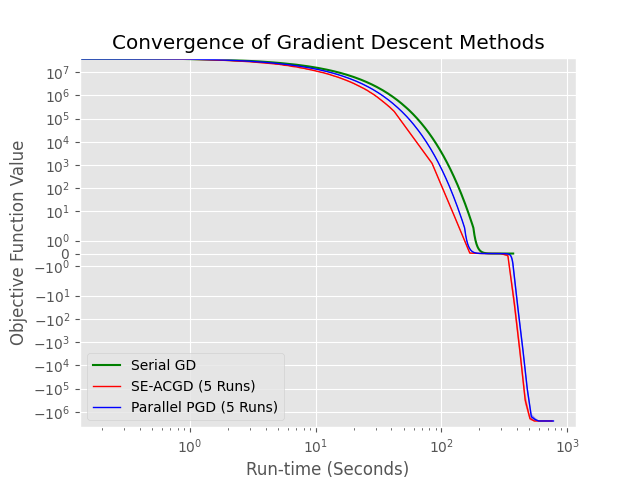}
         \caption{1e7 Dimensions, 8 Workers, No Delay.}
         \label{fig:8w-10M}
     \end{subfigure}
        \caption{Escaping Saddle Points Across Increasing Dimensionality.}
        \label{fig:escaping-saddle-comp}
\end{figure}

As one can see in Figure \ref{fig:escaping-saddle-comp}, SE-ACGD reaches the local optimum value of $-d/4$ for all dimensions.
In comparison, regular gradient descent (the green line) gets stuck at the saddle point with an objective value of 0.
Like SE-ACGD, perturbed gradient descent also manages to escape the saddle point and reach the local optimum.

Compared to parallel perturbed gradient descent, SE-ACGD manages to reach the local optimum \textit{faster even with no delay}.
This is an impressive feat, as comparing the two perturbed algorithms under no delay is to SE-ACGD's disadvantage: SE-ACGD performs best when there are parallelization delays (shown in \ref{app:parallelization-delays}). 
The superiority of SE-ACGD remains even as the dimensionality increases, showcasing the effectiveness of SE-ACGD at solving large-scale optimization problems.

\subsection{Parallelization Delays}
\label{app:parallelization-delays}

In this experiment, we showcase how effective SE-ACGD performs under parallelization delays compared to synchronous algorithms.
We perform this experiment at a large dimension, $d$ = 1 million, to simulate a large-scale problem.
To simulate an environment with parallelization delays, we inject an artificial delay at each iteration.
This artificial parallelization delay is added to one worker for each iteration, forcing the worker to sit idle during the duration of this delay.
\begin{figure}[!htbp]
     \centering
     \begin{subfigure}[b]{0.475\textwidth}
         \centering
         \includegraphics[width=0.75\textwidth]{Figures/8-Worker/8W-1000000-Dimension-Objective-0.000000.png}
         \caption{No Delay, 1e6 Dimensions, 8 Workers.}
         \label{fig:8w-nodelay}
     \end{subfigure}
     \hfill
     \begin{subfigure}[b]{0.475\textwidth}
         \centering
         \includegraphics[width=0.75\textwidth]{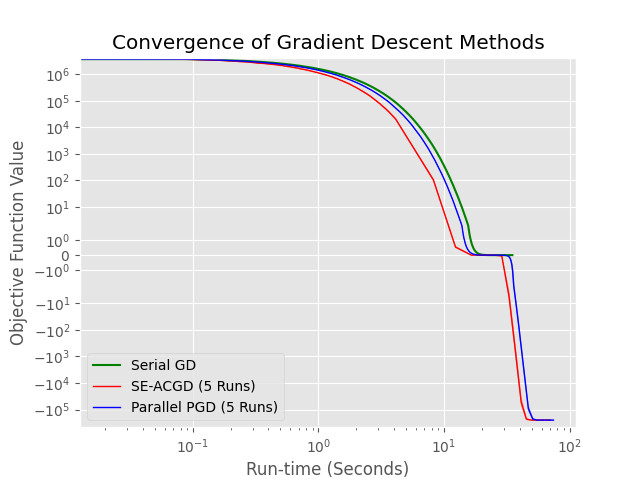}
         \caption{1e-2 Exp. Delay, 1e6 Dimensions, 8 Workers.}
         \label{fig:8w-0.001delay}
     \end{subfigure}
     \hfill
     \begin{subfigure}[b]{0.475\textwidth}
         \centering
         \includegraphics[width=0.75\textwidth]{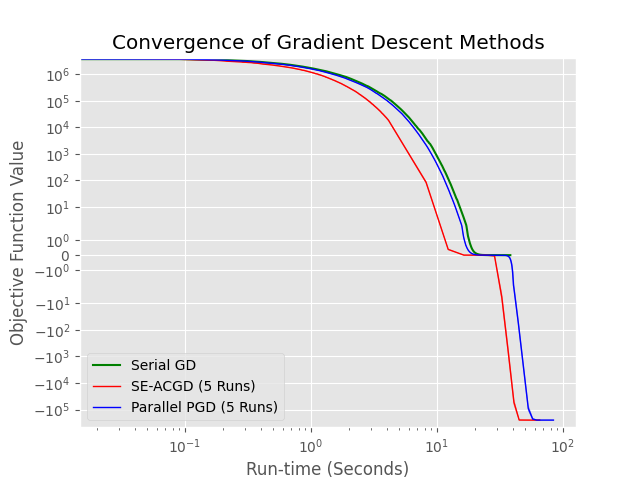}
         \caption{1e-1 Exp. Delay, 1e6 Dimensions, 8 Workers.}
         \label{fig:8w-0.01delay}
     \end{subfigure}
     \hfill
     \begin{subfigure}[b]{0.475\textwidth}
         \centering
         \includegraphics[width=0.75\textwidth]{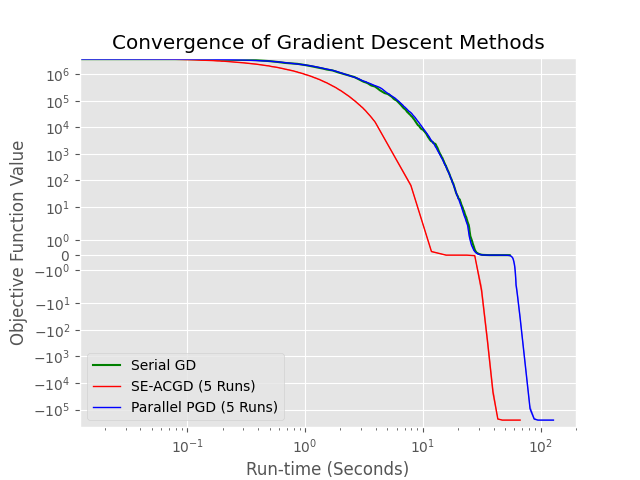}
         \caption{5e-2 Exp. Delay, 1e6 Dimensions, 8 Workers.}
         \label{fig:8w-0.05delay}
     \end{subfigure}
        \caption{Convergence Time Under Increasing Parallelization Delays.}
        \label{fig:delay-comp}
\end{figure}
We generate the delay randomly each iteration from an exponential distribution with the rate parameter $\lambda$ equal to the inverse of an inputted expected delay.
Within this experiment, we show how convergence time is effected as the expected delay is increased.

One can see in Figure \ref{fig:delay-comp} that the convergence time does indeed increase as the expected delay is increased.
However, this increase is much more dramatic for the synchronous algorithms than SE-ACGD.
While SE-ACGD and parallel perturbed gradient descent converge in a similar manner under no delay (Figure \ref{fig:8w-nodelay}), SE-ACGD leaps ahead of parallel perturbed gradient descent when the expected delay becomes 5e-2 (Figure \ref{fig:8w-0.05delay}).

These experimental results backup our claim about the effectiveness of SE-ACGD when experiencing parallelization delays.
Since SE-ACGD is already \textit{faster} than other synchronous algorithms under no delay (see \ref{app:saddle-evasion-exp}), and under delays it is no surprise that it vastly outperforms synchronous algorithms.
The results showcase that SE-ACGD is a more robust algorithm than its synchronous counterparts: it outperforms under ideal conditions and widens the gap under adverse conditions.

\subsection{Scalability}

Our final experiment is to show how SE-ACGD scales when varying the number of workers.
We again perform a large-scale test, with $d$ = 1 million.
Below we showcase how all parallel algorithms converge when using differing numbers of workers.

\begin{figure}[!htbp]
     \centering
     \begin{subfigure}[b]{0.475\textwidth}
         \centering
         \includegraphics[width=0.75\textwidth]{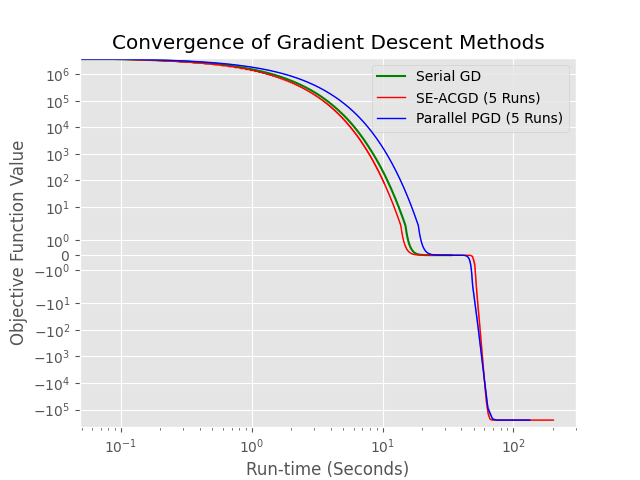}
         \caption{2 Workers, 1e6 Dimensions, No Delay.}
         \label{fig:2worker}
     \end{subfigure}
     \hfill
     \begin{subfigure}[b]{0.475\textwidth}
         \centering
         \includegraphics[width=0.75\textwidth]{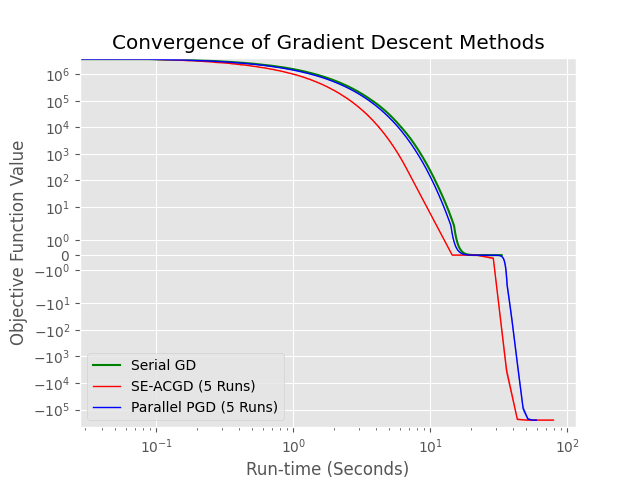}
         \caption{4 Workers, 1e6 Dimensions, No Delay.}
         \label{fig:4worker}
     \end{subfigure}
     \hfill
     \begin{subfigure}[b]{0.475\textwidth}
         \centering
         \includegraphics[width=0.75\textwidth]{Figures/8-Worker/8W-1000000-Dimension-Objective-0.000000.png}
         \caption{6 Workers, 1e6 Dimensions, No Delay.}
         \label{fig:6worker}
     \end{subfigure}
     \hfill
     \begin{subfigure}[b]{0.475\textwidth}
         \centering
         \includegraphics[width=0.75\textwidth]{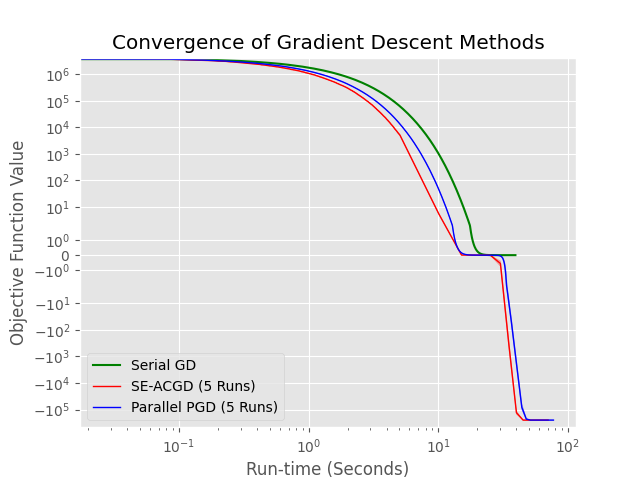}
         \caption{8 Workers, 1e6 Dimensions, No Delay.}
         \label{fig:8worker}
     \end{subfigure}
        \caption{Convergence Time With Varying Numbers of Workers.}
        \label{fig:vary-workers}
\end{figure}

As expected, the parallel algorithms begin to leave behind the serial gradient descent algorithm as the number of workers increase.
Thus, SE-ACGD performs better and scales with more workers.
In future work, we aim to scale up the number of workers even further on a cluster of CPUs. 
This will provide a greater reduction in the run-time of SE-ACGD on large-scale problems.

\subsection{Computational Devices} 

We ran all experiments locally on an Apple M1 Max Macbook Pro. 
SE-ACGD is built in Python and its workers communicates with one another via Open MPI, using MPI4Py (python bindings for MPI).
\vspace{-0.5em}
\section{Conclusion}
\label{sec:conclusion}
\vspace{-0.5em}
In this paper, we present a first-order asynchronous coordinate gradient descent algorithm which efficiently minimizes non-convex functions.
To the best of our knowledge, this is the first local optima convergence result for first-order asynchronous algorithm in non-convex settings.
We achieve this with only minor dependence on dimension, similar to \citep{jin2017escape, jin2017accelerated}.
Furthermore, we establish a relationship between bounded delays and convergence to second-order stationary points. 
Our work proves that asynchronous algorithms perform better than synchronous algorithms under large delays.

\section*{Ackowledgement}

Bornstein and Huang are supported by National Science Foundation NSF-IIS-FAI Award 2147276, Office of Naval Research, Defense Advanced Research Projects Agency Guaranteeing AI Robustness against Deception (GARD), and Adobe, Capital One and JP Morgan faculty fellowships.

\bibliography{bib}
\bibliographystyle{iclr2023_conference}

\newpage

\appendix

\newpage
\appendix
\onecolumn
\begin{center}
\LARGE{Appendix: \mytitlerunning}
\end{center}

\section{Parameters and Inequalities}
\label{appendix:parameters}

The user-inputted parameters are $\epsilon$, $\tau$, $L$, $\rho$, $\delta$, $d$, $x_0$, $W$, $\mu$ and $f^*$. Let $\Delta_f = f(x_0) - f^*$. Note that $d$, $W$, and $\mu$ are all greater than or equal to 1. The first set of hyper-parameters are defined as:
\begin{align}
    \label{eq:parameters1}
    \sigma = \max \{\frac{1280\sqrt{d} \Delta_f L \tau}{\sqrt{\pi}\epsilon^2 \delta} ,\; 8 \}, \quad \iota = \mu \cdot \text{log}_2 (\sigma), \quad \chi = \max \{ 1, \frac{\sqrt{\rho \epsilon}}{L^2} \}.
\end{align}
Within our algorithm, the learning rate $\eta$ is defined as
\begin{equation}
    \label{eq:learning_rate}
    \eta = \frac{1}{2L\tau^{1/2-\beta} \iota \chi}.
\end{equation}
The parameter $\beta > 0$ above is selected to be the largest possible value
\begin{equation}
\begin{aligned}
\label{eq:beta}
    &\max_\beta \quad \{\beta \; | \; \frac{15}{8}\tau^{1/2 - \beta} - \sqrt{\tau} - \frac{1}{2} \geq 0, \; \beta \leq 1/2\}.
\end{aligned}
\end{equation}
Finally, the remainder of our hyper-parameters are presented below
\begin{equation}
\begin{aligned}
\label{eq:parameters2}
& r = \eta \epsilon L, \quad \phi = \frac{5 \epsilon}{4 L\tau^{1/2-\beta} \iota \chi}, \quad \F = L(\frac{1}{\eta L} - \sqrt{\tau} - \frac{1}{2})\eta^2 \epsilon^2,\\
&\gamma = \frac{\delta \F}{\Delta_f},
\quad T = \frac{\text{log}_2(\sigma \iota^2 \chi^2)}{\eta \sqrt{\rho \epsilon}}, \quad r_0 = \frac{r \gamma \sqrt{\pi}}{2 \sqrt{d}},\quad M = \eta \epsilon.
\end{aligned}
\end{equation}
Before getting to the proofs, we mention that our analysis uses the following standard inequalities. For any $x_1, x_2, \ldots, x_N \in \R^N$ and $\epsilon > 0$, it holds that
\begin{align}
\label{eq:ineq}
\langle  x_1, x_2 \rangle \leq ||x_1||_2 \cdot ||x_2||_2 \leq \frac{\epsilon}{2} ||x_1||_2^2 + \frac{1}{2\epsilon}||x_2||_2^2.
\end{align}
We also use the following inequality
\begin{align}
\label{eq:log_eq_og}
    \text{log}_2(\sigma \iota^2 \chi^2) \leq \iota \chi \text{log}_2(\sigma) = \iota \chi (\frac{\iota}{\mu}) = \frac{\iota^2 \chi}{\mu}.
\end{align}

From $\eta \leq \frac{1}{2L\tau^{1/2-\beta} \iota \chi}$, Equation \eqref{eq:beta}, and $\iota \geq 3$, we see:
\begin{equation}
    \label{eq:beta_eq}
    (\frac{1}{\eta L} - \sqrt{\tau} - \frac{1}{2}) \geq 2\tau^{1/2 - \beta}\iota \chi - \sqrt{\tau} - \frac{1}{2} \geq \frac{1}{8}\tau^{1/2 - \beta}\iota\chi \geq \frac{3}{8}
\end{equation}
Finally, using all the parameters and inequalities above we find:
\begin{equation}
\label{eq:sigma_iota}
\begin{aligned}
     \frac{4 \phi}{\eta r_0} &= \frac{10\epsilon}{r_0} =  \frac{20\epsilon \sqrt{d}}{\sqrt{\pi}r \gamma} = \frac{20\sqrt{d} \Delta_f}{\sqrt{\pi}\eta L \delta \F} = \frac{20\sqrt{d} \Delta_f}{(\frac{1}{\eta L} - \sqrt{\tau} - \frac{1}{2})\sqrt{\pi}\eta^3 L^2 \epsilon^2 \delta}\\
     &= \frac{160\sqrt{d} \Delta_f L (\tau^{1/2-\beta})^3 \iota^3 \chi^3}{(\frac{1}{\eta L} - \sqrt{\tau} - \frac{1}{2})\sqrt{\pi}\epsilon^2 \delta} \leq \frac{1280\sqrt{d} \Delta_f L (\tau^{1/2-\beta})^3 \iota^3 \chi^3}{(\tau^{1/2-\beta}\iota\chi)\sqrt{\pi}\epsilon^2 \delta} \leq \sigma \iota^2 \chi^2
\end{aligned}
\end{equation}








\section{Detailed Proof of Guaranteed Descent for Hamiltonian}

Below we provide the full proof of Lemma \ref{lem:bounded_delay_decrease}, which proves that with a small enough step size the discrete Hamiltonian can become monotonically decreasing.

\begin{customlem}{1}[Guaranteed Discrete Hamiltonian Decrease]
Under Assumption \ref{assump:strict_saddle_lipschitz}, let the sequence
$\{x^j\}_{j \geq 0}$ be generated by ACGD with delays bounded by $\tau \geq 1$. There exists $\beta > 0$ such that choosing the step size
    $\eta \leq \frac{1}{2L\tau^{1/2 - \beta}\iota \chi}$ 
results in: 
    $E_j - E_{j+1} \geq L(\frac{1}{\eta L} - \sqrt{\tau} - \frac{1}{2})||x^{j+1} - x^{j}||_2^2 > 0$.
\end{customlem}

\begin{proof} (Lemma~\ref{lem:bounded_delay_decrease}):

The asynchronous coordinate gradient descent update rule for updating block $i$ and non-updating blocks $e$ is reprinted below:
\begin{equation}
\label{eq1:ASCD_update}
    x_{i}^{j+1}=x_{i}^{j}- \eta \nabla_{i}f(\hat{x}^{j}_i), \quad x_{e}^{j+1}=x_{e}^{j}
\end{equation}
Let $\Delta^j$ be defined as the difference between consecutive iterates and $\Delta_i^j$ as:
\begin{equation}
    \label{eq2:delta_diff}
    \Delta^j := x^{j+1} - x^j, \quad \Delta_i^j := \begin{cases} x^{j+1}_i - x^j_i, \text{ for block $i$} \\ 0, \text{ for non-updating blocks $e$}\end{cases}
\end{equation}
Therefore, by Equations \eqref{eq1:ASCD_update} and \eqref{eq2:delta_diff}, we get the following relationships:
\begin{equation}
    \label{eq3:delta_2}
    \Delta^j = -\eta \nabla f(\hat{x}^{j}_i), \quad \Delta_i^j = \begin{cases} -\eta \nabla_{i}f(\hat{x}^{j}_i), \text{ for block $i$}\\ 0, \text{ for non-updating blocks $e$} \end{cases}
\end{equation}
Using Equation \eqref{eq3:delta_2}, we see that:
\begin{equation}
\label{eq4:delta_3}
-\langle \Delta^j_i, \nabla_i f(\hat{x}^j_i) \rangle = -\langle \Delta^j, \nabla f(\hat{x}^j_i) \rangle = \frac{1}{\eta}||\Delta^j||_2^2 \rightarrow -\frac{1}{\eta}||\Delta^j||_2^2 - \langle \Delta^j, \nabla f(\hat{x}^j_i) \rangle = 0
\end{equation}

From the definition of the Lyapunov function (Definition \ref{def:hamiltonian_discrete}), we have:
\begin{equation}
\label{eq5:lyapunov_diff1}
        \resizebox{.9\hsize}{!}{$E_j - E_{j+1} = f(x^j) - f(x^{j+1}) + \frac{L}{2\sqrt{\tau}}\sum_{i=j - \tau}^{j-1}(i - (j - \tau) + 1)||\Delta^i||_2^2 - \frac{L}{2\sqrt{\tau}}\sum_{i=j - \tau + 1}^{j}(i - (j - \tau))||\Delta^i||_2^2$}
\end{equation}
\begin{equation}
        \label{eq5:lyapunov_diff2}
        E_j - E_{j+1} = f(x^j) - f(x^{j+1}) + \frac{L}{2\sqrt{\tau}}\sum_{i=j - \tau}^{j-1}||\Delta^i||_2^2 - \frac{L\sqrt{\tau}}{2}||\Delta^j||_2^2
\end{equation}
Since $f$ is $L$-smooth, we have:
\begin{equation}
\label{eq6:lsmooth}
f(x^{j+1}) \leq f(x^j) + \langle \nabla f(x^j), \Delta^j \rangle + \frac{L}{2}||\Delta^j||_2^2
\end{equation}
Using Equations \eqref{eq4:delta_3} and \eqref{eq6:lsmooth}, we get the following inequality:
\begin{align}
\label{eq7:eq7}
f(x^{j+1}) &\leq f(x^j) + \langle \nabla f(x^j), \Delta^j \rangle + \frac{L}{2}||\Delta^j||_2^2 - \frac{1}{\eta}||\Delta^j||_2^2 - \langle \Delta^j, \nabla f(\hat{x}^j_i) \rangle\\
&= f(x^j) + \langle \Delta^j, \nabla f(x^j) \rangle - \langle \Delta^j, \nabla f(\hat{x}^j_i) \rangle + (\frac{L}{2} - \frac{1}{\eta})||\Delta^j||_2^2 \\
&= f(x^j) + \langle \Delta^j, \nabla f(x^j)-\nabla f(\hat{x}^j_i) \rangle + (\frac{L}{2} - \frac{1}{\eta})||\Delta^j||_2^2
\end{align}
Manipulating the inequality above by using Equation \eqref{eq:ineq} (by substituting $L\sqrt{\tau}$ for $\epsilon$), we get the result:
\begin{equation}
\label{eq8:eq8}
f(x^{j+1}) \leq f(x^j) + \frac{L\sqrt{\tau}}{2} ||\Delta^j||_2^2 + \frac{1}{2L\sqrt{\tau}}||\nabla f(x^j)-\nabla f(\hat{x}^j_i)||_2^2 + (\frac{L}{2} - \frac{1}{\eta})||\Delta^j||_2^2
\end{equation}
Once again, by $L$-smooth properties, this yields:
\begin{align}
\label{eq9:eq9}
f(x^{j+1}) &\leq f(x^j) + \frac{L^2}{2L\sqrt{\tau}}||x^j - \hat{x}^j_i||_2^2 + (\frac{L}{2} + \frac{L\sqrt{\tau}}{2} - \frac{1}{\eta})||\Delta^j||_2^2\\
&= f(x^j) + \frac{L}{2\sqrt{\tau}}||x^j - \hat{x}^j_i||_2^2 + (\frac{L}{2} + \frac{L\sqrt{\tau}}{2} - \frac{1}{\eta})||\Delta^j||_2^2
\end{align}
From the definition of bounded delays and the triangle inequality, we know the following result:
\begin{equation}
    \label{triangle_inequality_delay}
    ||x^j - \hat{x}^j_i||_2 \leq \sum_{i = j -\tau}^{j-1}||x^{i+1} - x^i||_2 = \sum_{i = j -\tau}^{j-1}||\Delta^i||_2
\end{equation}
Using Equation \eqref{triangle_inequality_delay}, we find:
\begin{equation}
f(x^{j+1}) \leq f(x^j) + \frac{L}{2\sqrt{\tau}}\sum_{i = j - \tau}^{j-1}||\Delta^i||_2^2 + (\frac{L}{2} + \frac{L\sqrt{\tau}}{2} - \frac{1}{\eta})||\Delta^j||_2^2
\end{equation}
This is rearranged to yield:
\begin{equation}
\label{eq10:eq10}
f(x^j) - f(x^{j+1}) + \frac{L}{2\sqrt{\tau}}\sum_{i = j - \tau}^{j-1}||\Delta^i||_2^2 \geq  (\frac{1}{\eta} - \frac{L}{2} - \frac{L\sqrt{\tau}}{2} )||\Delta^j||_2^2
\end{equation}
Combining Equation \eqref{eq5:lyapunov_diff2} with Equation \eqref{eq10:eq10} leads to the desired result:
\begin{align}
E_j - E_{j+1} &\geq (\frac{1}{\eta} - \frac{L}{2} - \frac{L\sqrt{\tau}}{2})||\Delta^j||_2^2 - \frac{L \sqrt{\tau}}{2}||\Delta^j||_2^2\\
&= (\frac{1}{\eta} - \frac{L}{2} - L \sqrt{\tau})||\Delta^j||_2^2 = L(\frac{1}{\eta L} - \frac{1}{2} - \sqrt{\tau})||\Delta^j||_2^2
\end{align}
This is formally presented as:
\begin{equation}
    \label{eq11:eq11}
    E_j - E_{j+1} \geq L(\frac{1}{\eta L} - \sqrt{\tau} - \frac{1}{2})||x^{j+1} - x^j||_2^2
\end{equation}
The step size is selected as $\eta \leq \frac{1}{2L\tau^{1/2-\beta} \iota \chi}$ within this Lemma for a specific value $\beta > 0$. Applying the inequality for $\eta$ within Equation \eqref{eq11:eq11} results in:
\begin{equation}
\label{eq12:eq12}
E_j - E_{j+1} \geq L(2\tau^{1/2-\beta}\iota \chi - \sqrt{\tau} - \frac{1}{2})||x^{j+1} - x^j||_2^2
\end{equation}
Utilizing Equation \eqref{eq:beta_eq}, one can see that Equation \eqref{eq12:eq12} is always decreasing:
\begin{equation}
E_j - E_{j+1} \geq  \frac{3}{8}L||x^{j+1} - x^j||_2^2 > 0
\end{equation}
As one can see, if $\beta = 0$, then we have:
\begin{equation}
E_j - E_{j+1} > L(2\tau^{1/2} - \sqrt{\tau} - \frac{1}{2})||x^{j+1} - x^j||_2^2 \geq \frac{1}{2}L||x^{j+1} - x^j||_2^2 > 0
\end{equation}

\end{proof}

\section{Detailed Proof of Main Theorem}
\label{appendix:Main-Theorem-Proof}

In this section, we prove our main result piece by piece. We begin by proving that our algorithm sufficiently descends along large gradient and large momentum regions. Then we present the more complex proofs of sufficiently decreasing along saddle point regions. Once we have the proof of sufficient descent along all regions, we include the proof of escaping from a saddle point with a high probability. Once this is all demonstrated, we can finally prove our main theorem.

\subsection{Proof of Sufficient Descent Along Large-Gradient Region}
\label{appendix:large_grad}

Below we present the result of Theorem \ref{thm:large_gradient}. This theorem proves that our ACGD algorithm makes sufficient progress within the large-gradient region scenario.

\begin{customthm}{4}[\textbf{Large Gradient Scenario}]
Let function $f$ satisfy Assumption \ref{assump:strict_saddle_lipschitz}, $\eta$ satisfy Lemma \ref{lem:bounded_delay_decrease}, and let $\tau \geq 1$. If $||\nabla f(x^j)||_2 > \epsilon$, then by running Algorithm \ref{alg:Algo2} we have:
    $E_{j-\tau} - E_{j+1} > \F$.
\end{customthm}

\begin{proof} (of Theorem~\ref{thm:large_gradient}):

From the results of Lemma \ref{lem:bounded_delay_decrease}, we also know that choosing a step size $\eta \leq \frac{1}{2L\tau^{1/2 - \beta} \iota \chi}$ results in $E_j - E_{j+1} \geq (\frac{1}{\eta L} - \sqrt{\tau} - \frac{1}{2})L||x^j - x^{j+1} ||^2_2$. By L-smoothness of $f$, 
\begin{equation}
L||x^j - \hat{x}^j_i||_2 \geq || \nabla f(x^j) - \nabla f(\hat{x}^j_i)||_2 \geq ||\nabla f(x^j)||_2 - ||\nabla f(\hat{x}^j_i)||_2
\end{equation}
This results in:
\begin{equation}
\label{eq:thm5_1}
L||x^j - \hat{x}^j_i||_2 + ||\nabla f(\hat{x}^j_i)||_2 \geq ||\nabla f(x^j)||_2 
\end{equation}
From the update rule in the paper (Equation \eqref{eq:ACGD}), we know that $x^{j+1}_i = x^j_i - \eta \nabla_{i}f(\hat{x}^j_i)$. Rearranging this equation and taking the norm yields:
\begin{align}
    \label{eq:thm5_2}
     ||\nabla f(\hat{x}^j_i)||_2 = \frac{1}{\eta}||x^{j} - x^{j+1}||_2
\end{align}
Finally, from the delay of the gradient information and the triangle equation (Equation \eqref{triangle_inequality_delay}):
\begin{equation}
    \label{triangle_inequality_delay2}
    \sum_{k = j - \tau}^{j-1} ||x^k - x^{k+1}||_2 \geq ||\hat{x}^j_i - x^j||_2
\end{equation}
One can expand the sum over the delays as:
\begin{equation}
\label{eq:delay_sum}
\sum_{k = j - \tau}^j ||x^k - x^{k+1}||_2 = \sum_{k = j - \tau}^{j-1} ||x^k - x^{k+1}||_2 + ||x^j - x^{j+1}||_2
\end{equation}
From Equations \eqref{eq:thm5_2}, \eqref{triangle_inequality_delay2}, and \eqref{eq:delay_sum} (and multiplying through by $\frac{L}{\eta L}$) we create the following inequality:
\begin{align}
    \resizebox{.9\hsize}{!}{$\frac{L}{\eta L}\sum_{k = j - \tau}^{j} ||x^k - x^{k+1}||_2  \geq \frac{L||\hat{x}^j_i - x^{j}||_2}{\eta L} + \frac{L}{\eta L}||x^j - x^{j+1}||_2 \geq \frac{L||\hat{x}^j_i - x^{j}||_2}{\eta L} + ||\nabla f(\hat{x}^j_i)||_2$}
\end{align}
From defining $\eta \leq \frac{1}{2L\tau^{1/2 - \beta} \iota \chi}$, we see that $\eta L \leq \frac{1}{2\tau^{1/2 - \beta} \iota \chi}$. As $\tau, \iota,\chi \geq 1$, this leads to $\eta L \leq \frac{1}{2} < 1$ (as from the definition of $\beta$, 1/2 - $\beta \geq$ 0 for $\tau > 1$). 
From this result, the inequality above can be simplified to:
\begin{equation}
\frac{L||\hat{x}^j_i - x^{j}||_2}{\eta L} + ||\nabla f(\hat{x}^j_i)||_2 \geq L||\hat{x}^j_i - x^{j}||_2 + ||\nabla f(\hat{x}^j_i)||_2
\end{equation}
Now, using Equation \eqref{eq:thm5_1}, the inequality above transforms into:
\begin{align}
    &\frac{||\hat{x}^j_i - x^{j}||_2}{\eta} + ||\nabla f(\hat{x}^j_i)||_2 \geq L||\hat{x}^j_i - x^{j}||_2 + ||\nabla f(\hat{x}^j_i)||_2 \geq ||\nabla f(x^j)||_2
\end{align}
Finally we have reached the desired inequality:
\begin{equation}
    \sum_{k = j - \tau}^{j} ||x^k - x^{k+1}||_2 \geq \eta||\nabla f(x^j)||_2
\end{equation}
Squaring both sides leads to:
\begin{equation}
    (\sum_{k = j - \tau}^{j} ||x^k - x^{k+1}||_2)^2 \geq \eta^2 ||\nabla f(x^j)||_2^2
\end{equation}
Using Cauchy-Schwarz Inequality we have:
\begin{equation}
    \label{eq:cauchy_schwarz}
    \sum_{k = j - \tau}^{j} ||x^k - x^{k+1}||_2^2 \geq (\sum_{k = j - \tau}^{j} ||x^k- x^{k+1}||_2 )^2 \geq \eta^2 ||\nabla f(x^j)||_2^2
\end{equation}
Using the initial assumption that $||\nabla f(x^j)||_2 > \epsilon$ we find the desired inequality:
\begin{equation}
\label{eq:cauchy_worker}
    \sum_{k = j - \tau}^{j} ||x^k - x^{k+1}||_2^2 >  \eta^2 \epsilon^2
\end{equation}
Now, to finish the proof, we will utilize the assumptions from Lemma \ref{lem:bounded_delay_decrease}. For $k = j - \tau, \ldots, j$ we obtain that:
\begin{equation}
\sum_{k = j - \tau}^j E_{k} - E_{k+1} \geq \sum_{k = j - \tau}^j (\frac{1}{\eta L} - \sqrt{\tau} - \frac{1}{2})L||x^k - x^{k+1} ||^2_2
\end{equation}
Moving the constants outside of the sum yields:
\begin{equation}
\label{eq:thm5_3}
\sum_{k = j - \tau}^j E_{k} - E_{k+1}  \geq (\frac{1}{\eta L} - \sqrt{\tau} - \frac{1}{2})L \sum_{k = j - \tau}^j||x^k - x^{k+1} ||^2_2
\end{equation}
Equations \eqref{eq:cauchy_worker} and \eqref{eq:thm5_3} are used to get to the following inequality:
\begin{equation}
\sum_{k = j - \tau}^j E_{k} - E_{k+1}  > (\frac{1}{\eta L} - \sqrt{\tau} - \frac{1}{2})L \eta^2 \epsilon^2
\end{equation}
This is simplified down to become:
\begin{equation}
E_{j-\tau} - E_{j+1} > (\frac{1}{\eta L} - \sqrt{\tau} - \frac{1}{2})L \eta^2 \epsilon^2 = \F
\end{equation}

\end{proof}



\subsection{Proof of Sufficient Descent Along the Saddle-Point Region}

The direction of the minimum eigenvalue at a strict saddle points towards the descent direction, away from the saddle point. At some length along this direction, an iterate can escape from a saddle point. Our goal is to bound this length. 

We begin by proving what was stated informally above: at some length along the direction of the minimum eigenvalue an iterate can escape from a saddle point. This is summarized in the following Theorem.

\begin{lemma}
\label{lem:two_escape_points}
Under the conditions of \thm{saddle_point} , let $e_1$ be the smallest eigendirection of $\nabla^2 f(x^j)$. Consider two sequences $\{u^t\}_{t=0}^T$ and $\{w^t\}_{t=0}^T$, with $T \geq \frac{\text{log}_2(\sigma \iota^2 \chi^2)}{\eta \sqrt{\rho \epsilon}}$. Let these sequences be the iterates of ACGD starting from $u^0$ and $w^0$. We define $u^0=\tilde{x}+\xi$ and $w^0=u^0+\eta r_0 e_1$ ($\tilde{x}$ is the saddle point). The perturbation, $\xi$, comes from a uniform distribution over the ball with radius $\eta r$ centered at $\tilde{x}$. Then for $r_0 \in [\frac{r \gamma \sqrt{\pi}}{2 \sqrt{d}}, r]$, we have:
$\min\{E(u^T)-E(u^0), E(w^T)-E(w^0)\} < -4\F$.
\end{lemma}

Intuitively, \lem{two_escape_points} claims the following. Let $u^0$ be a point perturbed from a saddle point $x_s$ ($u^0$ is equivalent to $x_p$ in the proof sketch at the end of Section \ref{sec:Main_Result}). Now define the point $w^0$ to be a specified distance away from $u^0$ in the direction of the minimum eigenvalue ($w^0=u^0+\eta r_0 e_1$). 

Let each starting point, $u^0$ and $w^0$, undergo ACGD for a specified number of iterations $T$, generating sequences $\{u^t\}_{t=0}^T$ and $\{w^t\}_{t=0}^T$. \lem{two_escape_points} states that one of these two sequences will sufficiently decrease the Hamiltonian, and thus escape the saddle point.

\vspace{2em}
\textbf{Proof Sketch.} In order to prove \lem{two_escape_points}, we state the Localized-or-Improved property:

    \indent \textbf{1.} (Localized) We first show if ACGD from $u^0$ does not decrease the energy function enough, then all iterates must lie within a small ball around $u^0$.
    
    \indent \textbf{2.} (Improved) If ACGD starting from a point $u^0$ is stuck in a small ball around a saddle point, then ACGD from $w^0$ will decrease the energy function enough.
    
The Localized-or-Improved property is derived through the following two Lemmas. 

\begin{lemma}[\textbf{Localized}]\label{lem:localized}
Under the conditions of \lem{two_escape_points}, if $E(u^T)-E(u^0) \geq -4\F$, then for all $1\le t\le T$:
$\|u^t-\tilde{x}\|_2 \leq \phi$.
\end{lemma}

\begin{lemma}[\textbf{Improved}]\label{lem:improved}
Under the conditions of \lem{two_escape_points}, there exists a constant $\kappa > 0$ such that, if for all $1\le t\le T$, $\|u^t-\tilde{x}\|_2 \leq \phi$, then:
$E(w^T)-E(w^0) < -4\F$.
\end{lemma}

The first, Lemma \ref{lem:localized}, maintains that any instance of the objective function not decreasing by a sufficient amount implies that the iterates of ACGD all fall within a ball around the saddle point. This is the Localized portion of the property. 

The second, Lemma \ref{lem:improved}, states that if one sequence is stuck in a small ball around the saddle point, changing its starting point by $\eta r_0$ along the minimum eigenvalue direction will sufficiently decrease the objective function. These two Lemmas are the important pieces needed to prove \lem{two_escape_points}, and in turn finally prove \thm{saddle_point}.

\subsubsection{Proof of Localized-or-Improved Property}

The Localized-or-Improved property is proved below in Lemma \ref{lem:localized} (Localized) and Lemma \ref{lem:improved} (Improved). After proving both of these Lemmas, they can subsequently be used to prove  \lem{two_escape_points} that a specified perturbation at a saddle point will sufficiently decrease the Hamiltonian with high probability.

\begin{customlem}{7}[\textbf{Localized}]
Under the conditions of \lem{two_escape_points}, if $E(u^T)-E(u^0) \geq -4\F$, then for all $1\le t\le T$:
$\|u^t-\tilde{x}\|_2 \leq \phi$.
\end{customlem}

\begin{proof} (of Lemma~\ref{lem:localized}):

We are given that $E(u^T) - E(u^0) \geq - 4\F$, which is equivalent to $4\F \geq E(u^0) - E(u^T)$. Due to the monotonicity of $E$, $1 \leq t\leq T$, the triangle inequality, and $\eta \leq \frac{1}{2L\tau^{1/2 - \beta} \iota \chi}$, we have:
\begin{align}
4\F &\geq E(u^0) - E(u^T) \geq \sum_{j=0}^{t-1} E(u^j) - E(u^{j+1}) \\
&\geq \sum_{j=0}^{t-1} L(\frac{1}{\eta L} - \sqrt{\tau} - \frac{1}{2})||u^j - u^{j+1}||_2^2 \geq L(\frac{1}{\eta L} - \sqrt{\tau} - \frac{1}{2})||u^t - u^{0}||_2^2
\end{align}
Rearranging the inequality above yields:
\begin{align} \sqrt{\frac{4\F}{L(\frac{1}{\eta L} - \sqrt{\tau} - \frac{1}{2})}} \geq ||u^t - u^{0}||_2
\end{align}
By definition of the perturbation ball, we have $\eta r \geq ||u^0 - \tilde{x}||_2$. 
Now, we see that:
\begin{equation}
\resizebox{.9\hsize}{!}{$||u^t - \tilde{x}||_2 = ||u^t - u^0 + u^0 - \tilde{x}||_2 \leq ||u^t - u^0||_2 + ||u^0 - \tilde{x}||_2 \leq \sqrt{\frac{4\F}{L(\frac{1}{\eta L} - \sqrt{\tau} - \frac{1}{2})}} + \eta r$}
\end{equation}
Using the definitions of $\F$, $r$, $\eta$, and $\phi$ gets to the desired result:
\begin{align}
\sqrt{\frac{4\F}{L(\frac{1}{\eta L} - \sqrt{\tau} - \frac{1}{2})}} + \eta r = 2\sqrt{\eta^2 \epsilon^2} + \eta r = 2\eta \epsilon + \eta r = \eta(2\epsilon + r)\\
\leq \frac{1}{2L\tau^{1/2 - \beta}\iota \chi}(2\epsilon + r) = \frac{\epsilon}{2L\tau^{1/2 - \beta}\iota \chi}(2 + \eta L) \leq \frac{5\epsilon}{4L\tau^{1/2 - \beta} \iota \chi} = \phi
\end{align}
The final inequality above is found as $\eta L \leq \frac{1}{2}$. Thus the desired result is derived:
\begin{equation}
||u^t - \tilde{x}||_2 \leq \phi  
\end{equation}

\end{proof}

\begin{customlem}{8}[\textbf{Improved}]
Under the conditions of \lem{two_escape_points}, there exists a constant $\kappa > 0$ such that, if for all $1\le t\le T$, $\|u^t-\tilde{x}\|_2 \leq \phi$, then:
$E(w^T)-E(w^0) < -4\F$.
\end{customlem}

\begin{proof} (Lemma~\ref{lem:improved}):

\underline{Proof by Contradiction}

Without loss of generality, let $\tilde{x} = 0$ be the origin. The Hessian $H$ is then defined as
\begin{equation}
H = \nabla^2 f(\tilde{x}) = \nabla^2 f(0)
\end{equation}

In Lemma \ref{lem:localized}, it is assumed that $4\F \geq E(u^0) - E(u^T)$. Now, we will assume that there are two sequences, $\{u^t \}^T_{t = 0}$ and $\{w^t \}^T_{t = 0}$ (iterates of ACD from $u^0$ and $w^0$), that by Lemma \ref{lem:localized} for all $1\leq t \leq T$ both satisfy:
\begin{equation}
||u^t - \tilde{x}||_2 \leq \phi, \; ||w^t - \tilde{x}||_2 \leq \phi
\end{equation}
By the triangle equation, Equation \eqref{triangle_inequality_delay}, we have:
\begin{equation}
||\hat{u}^t_i - u^t||_2 \leq \sum_{k = t - \tau}^{t-1}||u^k - u^{k+1}||_2 \leq \sqrt{ \sum_{k = t - \tau}^{t-1}||u^k - u^{k+1}||_2^2}
\end{equation}
By Lemma \ref{lem:bounded_delay_decrease}, as the step size $\eta \leq \frac{1}{2L\tau^{1/2 - \beta} \iota \chi}$, we find that
\begin{align}
    \sqrt{ \sum_{k = t - \tau}^{t-1}||u^k - u^{k+1}||_2^2} &\leq \sqrt{\frac{1}{L(\frac{1}{\eta L} - \sqrt{\tau} - \frac{1}{2})}\sum_{k = t - \tau}^{t-1}(E(u^k) - E(u^{k+1}))}\leq \sqrt{\frac{E(u^0) - E(u^T)}{L(\frac{1}{\eta L} - \sqrt{\tau} - \frac{1}{2})}}
\end{align}
 
From the initial assumption ($4\F > E(x^0) - E(x^T)$), the definition of $\F$, and the definition of $\phi$, this leads to:
\begin{equation}
\sqrt{\frac{E(u^0) - E(u^T)}{L(\frac{1}{\eta L} - \sqrt{\tau} - \frac{1}{2})}} \leq \sqrt{\frac{4\F}{L(\frac{1}{\eta L} - \sqrt{\tau} - \frac{1}{2})}} = 2\sqrt{\epsilon^2 \eta^2} \leq  \frac{\epsilon}{L\tau^{1/2 - \beta} \iota \chi} \leq \frac{4}{5}\phi
\end{equation}
Therefore, we have the following result:
\begin{equation}
    ||\hat{u}^t_i - u^t||_2 \leq \frac{4}{5}\phi, \;  ||\hat{w}^t_i - w^t||_2 \leq \frac{4}{5}\phi
\end{equation}
Define $v^t = w^t - u^t$ and $\hat{v}^t_i = \hat{w}^t_i - \hat{u}^t_i$. Also, let $v^0 = \eta r_0 e_1$ as the initial distance between $u$ and $w$. Now, we can bound the following values:
\begin{align}
    ||v^t||_2 &= ||w^t - u^t||_2 = ||(w^t - \tilde{x}) - (u^t - \tilde{x})||_2\leq ||w^t - \tilde{x}||_2 + ||u^t - \tilde{x}||_2 \leq 2\phi
\end{align}
\begin{align}
    ||\hat{v}^t_i - v^t||_2 &= ||\hat{w}^t_i - \hat{u}^t_i - w^t + u^t||_2 = ||(\hat{w}^t_i - w^t) - (\hat{u}^t_i  - u^t)||_2\\
    &\leq ||\hat{w}^t_i - w^t||_2 + ||\hat{u}^t_i  - u^t||_2 \leq \frac{4}{5}\phi + \frac{4}{5}\phi = \frac{8}{5}\phi
\end{align}
Seeing that both equations are bounded above by $\phi$, we can find a value $\kappa > 0$ such that: 
\begin{equation}
||\hat{v}^t_i - v^t||_2 \leq \kappa||v^t||_2
\end{equation}
Consider the update equation for $w^t$:
\begin{align}
    u^{t+1} &+ v^{t+1} = w^{t+1} = w^t - \eta \nabla_i f (\hat{w}^t_i) = u^t + v^t - \eta \nabla_i f (\hat{u}^t_i + \hat{v}^t_i)\\
        &= u^{t+1} + \eta \nabla_i f(\hat{u}^t_i) + v^{t} - \eta \nabla_i f(\hat{u}^t_i + \hat{v}^t_i)\\
        v^{t+1} &= \eta \nabla_i f(\hat{u}_i^t) + v^{t} - \eta \nabla_i f(\hat{u}^t_i + \hat{v}^t_i)\\
        &= v^{t} + \eta \nabla_i f(\hat{u}_i^t) - \eta \nabla_i f(\hat{u}^t_i + \hat{v}^t_i) + \eta \nabla_i f(\hat{u}^t_i + v^t) - \eta \nabla_i f(\hat{u}_i^t + v^t)\\
        v^{t+1} &= v^{t} - \eta [(\int_0^1 \nabla \nabla_i f(\hat{u}^t_i + \theta v^t)d\theta)v^t] - \eta[\nabla_i f(\hat{u}^t_i + \hat{v}^t_i) - \nabla_i f(\hat{u}^t_i + v^t)]\\
        &= v^t - \eta [(H + \Delta_t)v^t] - \eta \zeta_t(\hat{v}^t_i - v^t)\\
         &= (I - \eta H - \eta \Delta_t)v^t - \eta \zeta_t(\hat{v}^t_i - v^t)\\
         &= (I - \eta H)v^t - \eta \sigma_t
\end{align}
The terms $\Delta_t$, $\zeta_t$, and $\sigma_t$ are defined as:
\begin{equation}
\label{eq:delta_t}
\Delta_t = \int_0^1 \nabla \nabla_i f(\hat{u}^t_i + \theta v^t)d\theta - H
\end{equation}
\begin{equation}
\label{eq:nabla_t}
\zeta_t = \int_0^1 \nabla \nabla_i f(\hat{u}^t_i + v^t + \theta(\hat{v}^t_i - v^t))d\theta
\end{equation}
\begin{equation}
\label{eq:sigma_t}
\sigma_t = -\eta[\Delta_t v^t + \zeta_t(\hat{v}^t_i - v^t)]
\end{equation}
The update equation for $v^{t+1}$ can now be rewritten using these terms as:
\begin{equation}
v^{t+1} = (I - \eta H)^{t+1}v^0 - \eta \sum_{m=0}^t (I - \eta H)^{t-m} \sigma_m
\end{equation}
Define the following terms $p(t+1)$ and $q(t+1)$ to get:
\begin{align}
&p(t+1) = (I - \eta H)^{t+1}v^0, \quad q(t+1) = \eta \sum_{m=0}^t (I - \eta H)^{t-m} \sigma_m\\
&v^{t+1} = p(t+1) - q(t+1)
\end{align}
By Hessian Lipschitz and Equations \eqref{eq:delta_t} and \eqref{eq:nabla_t} we have the following:
\begin{align}
||\Delta_t||_2 & = || \int_0^1 \nabla \nabla_i f(\hat{u}^t_i + \theta v^t)d\theta - H ||_2\\
&= || \int_0^1 \nabla^2 f(\hat{u}^t_i + \theta v^t)d\theta - \int_0^1 \nabla \nabla_e f(\hat{u}^t_i + \theta v^t)d\theta - H ||_2\\
&\leq \rho ||\hat{u}^t_i + v^t - \tilde{x}||_2 + ||\int_0^1 \nabla \nabla_e f(\hat{u}^t_i + \theta v^t)d\theta||_2 \\
&\leq \rho ||\hat{u}^t_i + v^t - \tilde{x}||_2 + ||\int_0^1 \nabla^2 f(\hat{u}^t_i + \theta v^t)d\theta||_2 \\
&\leq \rho ||\hat{u}^t_i + v^t - \tilde{x}||_2 + \rho || \hat{u}^t_i + v^t ||_2 = 2\rho || \hat{u}^t_i + v^t -\tilde{x} ||_2\\
&\leq 2\rho \big( || \hat{u}^t_i -\tilde{x} ||_2 + || v^t  ||_2 \big) \leq 2\rho \big( || \hat{u}^t_i - u^t||_2 + || u^t - \tilde{x} ||_2 + || v^t  ||_2 \big) \\ 
& \leq 2\rho \big( \frac{4}{5}\phi + \phi + 2\phi \big) = \frac{38}{5}\rho \phi
\end{align}
\begin{align}
||\zeta_t||_2 &= ||\int_0^1 \nabla \nabla_i f(\hat{u}^t_i + v^t + \theta(\hat{v}^t_i - v^t))d\theta||_2 \\
& \leq ||\int_0^1 \nabla^2 f(\hat{u}^t_i + v^t + \theta(\hat{v}^t_i - v^t))d\theta - \int_0^1 \nabla \nabla_e f(\hat{u}^t_i + v^t + \theta(\hat{v}^t_i - v^t))d\theta ||_2 \\
& \leq 2||\int_0^1 \nabla^2 f(\hat{u}^t_i + v^t + \theta(\hat{v}^t_i - v^t))d\theta||_2 \leq 2\rho ||\hat{u}^t_i + \hat{v}^t_i||_2 = 2\rho||\hat{w}^t_i||_2\\
&= 2\rho||\hat{w}^t_i - w^t + w^t - \tilde{x} + \tilde{x}||_2 \leq 2\rho(||\hat{w}^t_i - w^t||_2 + ||w^t - \tilde{x}||_2 + ||\tilde{x}||_2)\\
&\leq 2\rho(\frac{4}{5}\phi + \phi + 0) =  \frac{18}{5}\rho \phi
\end{align}
Now use induction to show that the term $q(t)$ is always small compared to the leading term $p(t)$. That
is, we wish to show:
\begin{equation}
||q(t)||_2 \leq \frac{||p(t)||_2}{2}, \quad \forall t \in T
\end{equation}
For the base case $t = 0$ this is true as $||q(0)||_2 = 0 \leq \frac{1}{2}||v^0||_2 = \frac{1}{2}||p(0)||_2$.
Now suppose that the induction claim is true up until $t$ ($0 \leq m \leq t$). We now see the following:
\begin{align}
&||v^m||_2 = ||p(m) - q(m)||_2 \leq ||p(m)||_2 + ||q(m)||_2 \leq \frac{3}{2}||p(m)||_2 = \frac{3}{2}(1 + \eta \sqrt{\rho \epsilon})^m \eta r_0\\
&||\hat{v}^m_i  - v^m||_2 \leq \kappa||v^m||_2 \leq \frac{3\kappa}{2}(1 + \eta \sqrt{\rho \epsilon})^m \eta r_0
\end{align}
Now we prove for the case $t+1 \leq T$. We start by bounding the norm of $\sigma_i$ using Equation \eqref{eq:sigma_t} as:
\begin{align}
||\sigma_m||_2  &= ||-\eta[\Delta_m v^m + \zeta_m(\hat{v}^m_i - v^m)]||_2 \leq \eta ( ||\Delta_m v^m||_2 + ||\zeta_m (\hat{v}^m_i - v^m)||_2)\\
&\leq \eta ( ||\Delta_m ||_2 ||v^m||_2 + ||\zeta_m||_2 || (\hat{v}^m_i - v^m)||_2)\\
&\leq \eta(\frac{57}{5}\rho \phi(1 + \eta \sqrt{\rho \epsilon})^m \eta r_0 + \frac{27\kappa}{5}\rho \phi(1 + \eta \sqrt{\rho \epsilon})^m \eta r_0)\\
&= \frac{27}{5}(\frac{57}{27} + \kappa)\eta \rho \phi (1 + \eta \sqrt{\rho \epsilon})^m \eta r_0
\end{align}
The upper bound of $q(t+1)$ is determined using the norm of $\sigma_i$ as:
\begin{align}
||q(t+1)||_2 &= \eta ||\sum_{m=0}^t (I - \eta H)^{t-m} \sigma_m||_2 \leq \eta ||\sum_{m=0}^t (I - \eta H)^{t-m}||_2 ||\sigma_m||_2
\end{align}
\begin{align}
&\leq \frac{27}{5}(\frac{57}{27} + \kappa)\eta^2 \rho \phi (1 + \eta \sqrt{\rho \epsilon})^m \eta r_0 ||\sum_{m=0}^t (I - \eta H)^{t-m}||_2
\end{align}
\begin{align}
&\leq \frac{27}{5}(\frac{57}{27} + \kappa)\eta^2 \rho \phi T (1 + \eta \sqrt{\rho \epsilon})^t \eta r_0\\
& \leq \frac{27}{5}(\frac{57}{27} + \kappa)\eta^2 \rho \phi T ||p(t+1)||_2 \leq \frac{1}{2}||p(t+1)||_2
\end{align}
The last inequality above ($\frac{27}{5}(\frac{57}{27} + \kappa)\eta^2 \rho \phi T \leq \frac{1}{2}$) follows from the definition of $T$, $\phi$, $\eta$, and $\F$:
\begin{align}
\frac{27}{5}(\frac{57}{27} + \kappa)\eta^2 \rho \phi T &= \frac{57}{5}(\frac{1}{2} + \kappa)\rho \eta^2(\frac{5\epsilon}{4L\tau^{1/2 - \beta} \iota \chi})(\frac{\text{log}_2(\sigma \iota^2 \chi^2)}{\eta \sqrt{\rho \epsilon}}) \\
&= \frac{27}{4}(\frac{57}{27} + \kappa)\frac{\eta \sqrt{\rho \epsilon} }{L\tau^{1/2 - \beta} \iota \chi}(\text{log}_2(\sigma \iota^2 \chi^2))\\
&= \frac{27}{4}(\frac{57}{27} + \kappa) \frac{\sqrt{\rho \epsilon}}{2 L^2 \tau^{1 - 2\beta} \iota^2 \chi^2}(\text{log}_2(\sigma \iota^2 \chi^2))\\
\label{eq:log_ineq}
&\leq \frac{\sqrt{\rho \epsilon}(\frac{27}{4}(\frac{57}{27} + \kappa))}{2 L^2 \tau^{1 - 2\beta} \iota^2 \chi^2}(\iota \chi (\frac{\iota}{\mu})) = \frac{\sqrt{\rho \epsilon}(\frac{27}{4}(\frac{57}{27} + \kappa))}{2 L^2 \tau^{1 - 2\beta} \chi \mu}\\
&\leq \frac{\sqrt{\rho \epsilon}(\frac{27}{4}(\frac{57}{27} + \kappa))}{2 L^2 (\frac{\sqrt{\rho \epsilon}}{L^2}) \mu} = \frac{27(\frac{57}{27} + \kappa)}{8 \mu} \leq \frac{1}{2}
\end{align}
The inequality in Equation \eqref{eq:log_ineq} holds by Equation \eqref{eq:log_eq_og} in Appendix \ref{appendix:parameters}. The last line follows by picking a large enough constant $\mu$ for the inequality to hold. This finishes the inductive proof showing that $||q(t)||_2 \leq \frac{||p(t)||_2}{2}$, for all $t \in T$. Using this result, we take the norm $||v^T||_2$ and see the following:
\begin{align}
&||v^T||_2 = ||p(T) - q(T)||_2 \geq ||p(T)||_2 - ||q(T)||_2 \geq \frac{1}{2}||p(T)||_2 = \frac{1}{2}(1 + \eta \sqrt{\rho \epsilon})^T \eta r_0\\
&||v^T||_2 = \frac{1}{2}(1 + \eta \sqrt{\rho \epsilon})^T \eta r_0 > 2^{T \eta \sqrt{\rho \epsilon}}(\frac{\eta r_0}{2})
\end{align}
The last line follows from the inequality $(1+a)^x > 2^{ax}$ for $x > 0$ and $a \in (0,1)$. As $\eta \leq \frac{1}{2L\tau^{1/2 - \beta}\iota \chi}$, we see that $\eta \sqrt{\rho \epsilon} \leq \frac{\sqrt{\rho \epsilon}}{2L\tau^{1/2 - \beta}\iota \chi} \leq \frac{1}{2\tau^{1/2 - \beta}\iota \chi} \leq \frac{1}{2}$ and thus is contained within the interval (0,1). Since $T = \frac{\text{log}_2(\sigma \iota^2 \chi^2)}{\eta \sqrt{\rho \epsilon}}$, and using Equation \eqref{eq:sigma_iota}, we find the desired contradiction:
\begin{align}
||v^T||_2 > 2^{\eta \sqrt{\rho \epsilon} T} (\frac{\eta r_0}{2}) = (\sigma \iota^2 \chi^2)(\frac{\eta r_0}{2}) \geq (\frac{4\phi}{\eta r_0})(\frac{\eta r_0}{2})= 2\phi
\end{align}
\end{proof}

\subsubsection{Proof of Bounded Stuck Region}

Finally, using the proof of the Localized-or-Improved property, one can prove \lem{two_escape_points}:

\begin{customlem}{6}
Under the conditions of \thm{saddle_point} , let $e_1$ be the smallest eigendirection of $\nabla^2 f(x^j)$. Consider two sequences $\{u^t\}_{t=0}^T$ and $\{w^t\}_{t=0}^T$, with $T \geq \frac{\text{log}_2(\sigma \iota^2 \chi^2)}{\eta \sqrt{\rho \epsilon}}$. Let these sequences be the iterates of ACGD starting from $u^0$ and $w^0$. We define $u^0=\tilde{x}+\xi$ and $w^0=u^0+\eta r_0 e_1$ ($\tilde{x}$ is the saddle point). The perturbation, $\xi$, comes from a uniform distribution over the ball with radius $\eta r$ centered at $\tilde{x}$. Then for $r_0 \in [\frac{r \gamma \sqrt{\pi}}{2 \sqrt{d}}, r]$, we have:
$\min\{E(u^T)-E(u^0), E(w^T)-E(w^0)\} < -4\F$.
\end{customlem}

\begin{proof} (Lemma~\ref{lem:two_escape_points}):

Assume the non-trivial case in which the sequence $\{u^t\}_{t=0}^T$ (iterates of asynchronous coordinate gradient descent starting at $u^0$) is stuck at a saddle point, and thus $E(u^T) - E(u^0) \geq -4\F$. Then, by Lemma \ref{lem:improved}, we know that the sequence $\{w^t\}_{t=0}^T$ (iterates of asynchronous coordinate gradient descent starting at $w^0 = u^0 + \eta r_0 e_1$) will escape the saddle point. Therefore, we would have that $E(w^T) - E(w^0) < -4\F$. It directly follows that:
\begin{equation}
    \min\{E(u^T)-E(u^0), E(w^T)-E(w^0)\} < -4\F.
\end{equation}

\end{proof}

\subsubsection{Proof of Sufficient Descent At Saddle Point With High Probability}

Below, in Theorem \ref{thm:saddle_point}, we prove that with probability $1-\gamma$, our ACGD algorithm will sufficiently decrease the Hamiltonian when perturbing a point stuck at a saddle point. 

\begin{customthm}{5}[\textbf{Saddle Point Scenario}]
Let function $f$ satisfy Assumption \ref{assump:strict_saddle_lipschitz} and let $\eta$ satisfy Lemma \ref{lem:bounded_delay_decrease}. If $\|\nabla f(x^j)\|_2 \leq \epsilon$ and $\lambda_{\min}(\nabla^2 f(x^j)) < -\sqrt{\rho\epsilon}$, then $x^j$ is located at a saddle point. Let $y^0=x^j+\xi$, where $\xi$ comes from the uniform distribution over ball with radius $\eta r$, and $\{y^t\}_{t=0}^T$ are the iterates of ACGD starting from $y^0$ with $T \geq \frac{\log_2(\frac{4\phi}{\eta r_0})}{\eta \sqrt{\rho \epsilon}}$. Let $x^{j+1} = y^T$. Then, with at least probability $1-\gamma$, running Algorithm \ref{alg:Algo3} once results in:
    $E_{j} - E_{j+1} > \F$.
\end{customthm}

\begin{proof} (Theorem~\ref{thm:saddle_point}):

Denote $\tilde{x}$ as an iterate stuck at a saddle point. By selecting $\eta \leq \frac{1}{2L\tau^{1/2 - \beta}\iota \chi}$, we know from Equation \eqref{eq:beta_eq} that $(\frac{1}{\eta L} - \sqrt{\tau} - \frac{1}{2}) \geq \frac{3}{8}$ and thus $\F \geq \frac{3L\eta^2 \epsilon^2}{8} $. By applying Lemma \ref{lem:two_escape_points}, we know if $\tilde{x}$ is stuck at a saddle point then $\tilde{x}$ can be perturbed in a manner ($y^0 = \tilde{x} + \eta r_0 e_1$) such that the perturbed point $y^0$ is guaranteed not to be stuck ($E(y^T) - E(y^0) < -4\F$) after $T$ iterations of ACGD. At worst case, adding a perturbation $\xi$ will increase the function value ($L$-smooth properties, small gradient, definition of $r$, definition of $\eta$, and $\xi \le \eta r$) by:
\begin{align}
    f(y^0) \leq f(\tilde{x}) + \nabla f(\tilde{x})^T \xi + \frac{L||\xi||^2}{2} \leq f(\tilde{x}) + \epsilon \eta r + \frac{L||\eta r||^2}{2} \\
    f(y^0) - f(\tilde{x}) \leq \epsilon \eta r + \frac{L\eta^2 r^2}{2} = \epsilon^2 \eta^2 L + \frac{L^3 \eta^4 \epsilon^2}{2} < \epsilon^2 \eta^2 L + \frac{L \eta^2 \epsilon^2}{8}\\
    f(y^0) - f(\tilde{x}) <  \frac{8}{3}\F   + \frac{1}{3}\F = 3\F
\end{align}
Thus we have a desired result that:
\begin{equation}
    \label{eq:thm6_1}
    E(y^0) - E(\tilde{x}) = f(y^0) - f(\tilde{x}) < 3\F
\end{equation}

The $d$-dimensional perturbation ball centered at $\textbf{x}$ with radius $\eta r$ is denoted as $B_{\textbf{x}}^{(d)}(\eta r)$. In Lemma \ref{lem:improved} it is defined that $w^0 = u^0 + \eta r_0e_1$. Since $r_0 \geq \frac{r \gamma \sqrt{\pi}}{2\sqrt{d}}$, for any two points along the $e_1$ direction that are at least $\frac{r \gamma \sqrt{\pi}}{\sqrt{d}}$ away from each other, we now know one must not be in the stuck region $R_{stuck}.$ Therefore, the stuck region can be bounded. If there is a point $\tilde{u} \in R_{stuck}$ along $e_1$, then the stuck interval around $\tilde{u}$ must be of length $\frac{r \gamma \sqrt{\pi}}{\sqrt{d}}$. This is the thickness of $R_{stuck}$ along the $e_1$ direction. Using calculus, this can be turned into an upper bound on the volume of the stuck region.

As shown in \cite{jin2017escape}, denote $I_{stuck}(\cdot)$ to be the indicator function of being in $R_{stuck}$ and $x = (x^{(1)},\bar{x})$ where $x^{(1)}$ is in the $e_1$ direction and $\bar{x}$ is the remaining $d-1$ dimensional vector. The volume of $R_{stuck}$ is determined as:
\begin{align}
Vol(R_{stuck}) &= \int_{B^{(d)}_{\tilde{x}}(\eta r)}dx \cdot I_{stuck}(x)\\
&= \int_{B^{(d-1)}_{\tilde{x}}(\eta r)}d\bar{x} \int_{x^{(1)} - \sqrt{(\eta r)^2 - ||\bar{\tilde{x}} - \bar{x} ||^2}}^{x^{(1)} - \sqrt{(\eta r)^2 - ||\bar{\tilde{x}} - \bar{x} ||^2}} dx^{(1)} \cdot I_{stuck}(x)
\end{align}
Since we know that $r_0 \geq \frac{r \gamma \sqrt{\pi}}{2\sqrt{d}}$
\begin{align}
Vol(R_{stuck}) \leq \int_{B^{(d-1)}_{\tilde{x}}(\eta r)}d\bar{x} \cdot  \bigg(\frac{ 2\eta r \gamma \sqrt{\pi}}{2\sqrt{d}} \bigg) = Vol(B_{\tilde{x}}^{(d-1)}(\eta r))\bigg(\frac{\eta r \gamma \sqrt{\pi}}{\sqrt{d}}\bigg)
\end{align}
The volume ratio of the stuck region can now be computed as:
\begin{equation}
\frac{Vol(R_{stuck})}{Vol(B_{\tilde{x}}^{(d)}(\eta r))} \leq \frac{Vol(B_{\tilde{x}}^{(d-1)}(\eta r))\bigg(\frac{\eta r \gamma \sqrt{\pi}}{\sqrt{d}}\bigg)}{Vol(B_{\tilde{x}}^{(d)}(\eta r))}
\end{equation}
The volume of a $d$ dimensional ball of radius $r$ is given as
\begin{equation}
    V_d(r)=\frac{\pi^{\frac{d}{2}}}{\Gamma (\frac{d}{2} + 1)}r^d
\end{equation}
Therefore, the ratio of volumes can be simplified to
\begin{align}
\frac{Vol(R_{stuck})}{Vol(B_{\tilde{x}}^{(d)}(\eta r))} \leq \frac{\bigg(\frac{\eta r \gamma \sqrt{\pi}}{\sqrt{d}}\bigg)\frac{\pi^{\frac{d-1}{2}}}{\Gamma (\frac{d-1}{2} + 1)}(\eta r)^{d-1}}{\frac{\pi^{\frac{d}{2}}}{\Gamma (\frac{d}{2} + 1)}(\eta r)^d} \\ 
=\frac{\bigg(\frac{\eta r \gamma \sqrt{\pi}}{\sqrt{d}} \bigg) \Gamma (\frac{d}{2} + 1)}{\eta r\sqrt{\pi}\;\Gamma (\frac{d-1}{2} + 1)}
= \bigg( \frac{\gamma}{\sqrt{d}}\bigg) \frac{\Gamma (\frac{d}{2} + 1)}{\Gamma (\frac{d}{2} + 1/2)}
\end{align}
By property of the gamma function, $\frac{\Gamma (x + 1)}{\Gamma (x + 1/2)} \leq \sqrt{x + \frac{1}{2}}$. This is used to simplify the inequality above to become:
\begin{equation}
\frac{Vol(R_{stuck})}{Vol(B_{\tilde{x}}^{(d)}(\eta r))} \leq \bigg( \frac{\gamma}{\sqrt{d}}\bigg) \sqrt{\frac{d}{2} + \frac{1}{2}} \leq  \bigg( \frac{\gamma}{\sqrt{d}}\bigg) \sqrt{d} = \gamma
\end{equation}
Therefore we see that
\begin{equation}
    \frac{Vol(R_{stuck})}{Vol(B_{\tilde{x}}^{(d)}(\eta r))}\leq \gamma
\end{equation}
This result shows that if $\tilde{x}$ is originally stuck at a saddle point, then the added perturbation results in $y^0 \notin R_{stuck}$ with a probability of $1-\gamma$. In this case, we have (by Equation \eqref{eq:thm6_1} and Lemma \ref{lem:two_escape_points}):
\begin{align}
E_{j+1} - E_{j} = E(y^T) - E(\tilde{x})& = E(y^T) - E(y^0) + E(y^0) - E(\tilde{x}) < -4\F + 3\F = -\F
\end{align}
The saddle point $\tilde{x}$ is equivalent to the stuck iterate $x^j$, and the output of the perturbation algorithm $y^T$ is equivalent to $x^{j+1}$. This leads to our desired result, that with a probability of $1-\gamma$ we have:
\begin{equation}
    E_j - E_{j+1} > \F
\end{equation}

\end{proof}

\subsubsection{Proof of Main Theorem}

Finally, we present the proof of our main theorem below. This theorem states that with high probability $1-\delta$, our ACGD efficiently escapes from saddle points and converges to $\epsilon$-second-order stationary points. Furthermore, in our main theorem, Theorem \ref{thm:main_theorem}, we prove the convergence rate described in our work above.

\begin{customthm}{3}[\textbf{Main Theorem}]
Let the function $f$ satisfy Assumption \ref{assump:strict_saddle_lipschitz}. For any $\delta>0$ and $\epsilon\le\frac{L^2}{\rho}$, there exists a step size $\eta$ satisfying Lemma \ref{lem:bounded_delay_decrease} such that Algorithm \ref{alg:MainAlgo} will output an $\epsilon$-second-order stationary point, with probability $1-\delta$, in the following number of iterations:
$$
    O\bigg(\frac{L^2 \Delta_f \tau^{1-2\beta}}{\epsilon^{2.5}\rho^{0.5}} \text{log}_2^4(\frac{d \tau L \Delta_f}{\delta \epsilon})\bigg)
$$
\end{customthm}

\begin{proof} (of Theorem~\ref{thm:main_theorem}):

The maximum number of iterations of the algorithm (if every point is a saddle point) is defined as:
\begin{equation}
T_{max} =  \frac{T\Delta_f}{\F}
\end{equation}
This worst case scenario takes $\frac{T(f(x_0) - f^*)}{\F}$ iterations. Therefore, this is on the order of:
\begin{align}
\frac{T(f(x_0) - f^*)}{\F} &= \frac{\bigg( \frac{\text{log}_2( \sigma \iota^2 \chi^2)}{\eta \sqrt{\rho \epsilon}} \bigg)(f(x_0) - f^*)}{ \eta^2 \epsilon^2 L(\frac{1}{\eta L} - \sqrt{\tau} - \frac{1}{2})}\\
&=\frac{\text{log}_2( \sigma \iota^2 \chi^2) \Delta_f}{\eta^3 \epsilon^2 L (\frac{1}{\eta L} - \sqrt{\tau} - \frac{1}{2}) \sqrt{\rho \epsilon}} \leq  \frac{ 8(\iota^2 \chi) \Delta_f }{\mu \eta^3 \epsilon^2 L(\tau^{1/2-\beta}\iota \chi) \sqrt{\rho \epsilon}}\\
&= O\bigg(\frac{L^2 \tau^{1 - 2\beta} \Delta_f}{ \epsilon^2 \sqrt{\rho \epsilon}} \text{log}_2^4(\frac{d \tau  L \Delta_f}{\delta \epsilon})\bigg)\\
&= \tilde{O}\bigg(\frac{L^2 \Delta_f  \tau^{1-2\beta}}{\epsilon^{2.5}\rho^{0.5}}\bigg)
\end{align}

From Theorem \ref{thm:saddle_point}, the probability that the saddle point escapes during one run of the perturbation algorithm (with $T$ iterations) is $1 - \gamma$. In the worst case scenario described above, the perturbation algorithm would be run $\frac{\Delta_f}{\F}$ times. Therefore, the probability of failure for the entire algorithm is described as $(1 - \gamma)^{\frac{\Delta_f}{\F}}$. Since $0 \leq \gamma \leq 1$ (it is a probability) and $\frac{\Delta_f}{\F} \geq 1$ (there must be at least one run of Algorithm \ref{alg:MainAlgo}), one can see that the following inequality holds:
\begin{equation}
    (1 - \gamma)^{\frac{\Delta_f}{\F}} \geq 1 - \frac{\Delta_f}{\F}\gamma
\end{equation}
Now, to complete the proof, one must show that $1 - \frac{\Delta_f}{\F}\gamma \geq 1 - \delta$. Using the definition of $\gamma$, this simplifies to:
\begin{equation}
    \delta \geq \frac{\Delta_f}{\F}\gamma = \frac{\Delta_f}{\F}\frac{\delta \F}{\Delta_f} = \delta
\end{equation}
\end{proof}

\end{document}